\documentclass[review]{elsarticle}
\usepackage{amsmath}
\usepackage{amstext}
\usepackage{amssymb}
\usepackage{amsthm}
\usepackage{caption}
\usepackage{lineno,hyperref}
\usepackage{subfigure}
\usepackage{graphics}
\usepackage{booktabs}
\usepackage{color}
\modulolinenumbers[5]

\newtheorem{lemma}{Lemma}

\newtheorem{thm}{Theorem}
\newtheorem{cor}{Corollary}
\newtheorem{remark}{Remark}
\newtheorem{pro}{Proposition}

\begin{document}

\begin{frontmatter}

\title{A necessary and sufficient condition for bounds on the sum of a list of real numbers and its applications}


%
%
%
\author{Xiwu Yang}
\ead{cinema@lnnu.edu.cn}
\address{School of Mathematics, Liaoning Normal University, Dalian 116029,  P.R. China}
\begin{abstract}
Let $x_1,...,x_n$ be a list of real numbers, let $s :=\sum_{i=1}^{n}x_i$
and let $h:\mathbb{N} \rightarrow \mathbb{R}$ be a function.
We gave a necessary and sufficient condition for $s>h(n)$ (respectively, $s<h(n)$).
Let $G=(V,E)$ be a graph, let $\{H_1,...,H_n\}$ and $\{V_1,...,V_n\}$ be a decomposition and a partition of $G$, respectively.
Let $H_{i,j}$ and $V_{i,j}, i\leq j,$ be the union of $H_i,...,H_j$ and $V_i,...,V_j$, respectively, where subscripts are taken modulo $n$.
$G$ is \emph{generalized periodic} or \emph{partition-transitive} if for each pair of integers $(i,j)$,
$H_{i,i+k}$ and $H_{j,j+k}$ or $V_{i,i+k}$ and $V_{j,j+k}$ are isomorphic for all $k$, $1\leq k\leq n$, respectively.
Let $f:E \rightarrow \mathbb{R}$ and $g:V \rightarrow \mathbb{R}$ be mappings,
let the \emph{weight} of $f$ or $g$ on $G$ be $\Sigma_{e\in E}f(e)$ or $\Sigma_{v\in V}g(v)$, respectively.
Suppose that parameters $\lambda$ and $\xi$ of $G$ can be expressed as the minimum or maximum weight of specified $f$ and $g$, respectively.
Then our conditions imply a necessary and sufficient condition for $\lambda(G_1)=h(n)$ (respectively, $\xi(G_2)=h(n)$), where $G_1$ is generalized periodic and $G_2$ is partition-transitive.
For example, the crossing number $\textrm{cr}(\odot(T^n))$ of a periodic graph $\odot(T^n)$, $\textrm{cr}(\odot(T^n))=h(n)$.
As applications, we obtained $\textrm{cr}(C(4n;\{1,4\}))$ of the circulant $C(4n;\{1,4\})$,
the paired domination number of $C_5\Box C_n$
and the upper total domination number of $C_4\Box C_n$.
\end{abstract}
\end{frontmatter}


\section{Introduction}
Let $G=(V,E)$ be a graph,
let $f:E \rightarrow \mathbb{R}$ be a mapping
and let the \emph{weight} of $f$ on $G$ be $\Sigma_{e\in E}f(e)$.
For some terminology and notation not defined here, readers can refer to \cite{Bon}.
Let $\Lambda$ be the set of parameters that can be expressed as the minimum or maximum of the weight of a specified $f$ on $G$.
For example, the edge covering number, the matching number and both the crossing number and
the maximum crossing number in a given surface of a graph are members of $\Lambda$.
The \emph{crossing number} of $G$ in an arbitrary surface $\Sigma$, denoted by $\textrm{cr}_{\Sigma}(G)$, is the minimum number
of pairwise crossings of edges in a drawing of $G$ in $\Sigma$.
The \emph{maximum crossing number} \cite{Sch22} of $G$, max-$\textrm{cr}(G)$, is the largest number of crossings in any drawing of $G$
in which every pair of edges has at most one point in common (including a shared endpoint; touching points are forbidden.)

Robertson and Seymour \cite{Rob} showed that for each surface $\Sigma$ the class of graphs that embed into
$\Sigma$ can be characterized by a finite list Forb($\Sigma$) of minimal forbidden minors.
For the 2-sphere $\mathbb{S}_0$, Forb($\mathbb{S}_0$) is $\{K_5,K_{3,3}\}$ \cite{Kur}.
For the projective plane $\mathbb{N}_1$, Forb($\mathbb{N}_1$) is a list of 35 graphs \cite{Arc,Glo} and
$\mathbb{N}_1$ is the only other surface for which the complete list of forbidden minors is known \cite{Moh}.
Richter and \v{S}ir\'{a}\v{n} \cite{Ric} obtained the crossing number of $K_{3,n}$ in an arbitrary surface.
Shahrokhi et al. \cite{Sha} gave a general lower bound on the crossing number of a graph on a compact 2-manifold.
The crossing number of several graphs in torus \cite{Guy68,Guy69,Ho09}, in the projective plane \cite{Ris,Ho05} and in the Klein bottle \cite{Ris01} were obtained.

The known results on the maximum crossing numbers of graph classes in the plane is few.
It contains the complete graph \cite{Rin}, the complete $n$-partite graph \cite{Harb}, trees \cite{Pia} and cycle $C_n$ \cite{Woo71}.
In particular, let $\textrm{cr}(G)$ denote the crossing number of $G$ in the plane (or $\mathbb{S}_0$).
A graph is \emph{periodic} if it can be obtained by joining identical pieces in a circular fashion.
Let $\odot(T^n)$ be a periodic graph, where $T$ is a tile.
Pinontoan and Richter \cite{Pin03} showed that the limit
$c(T)={\lim_{n \to \infty}}\textrm{cr}(\odot(T^n))/n$ exists.
This result plays central role in constructions of crossing-critical graphs \cite{Pin04} and in explaining certain event.
Dvo\v{r}\'{a}k and Mohar \cite{RefDvo} proved that for every tile $T$ and for every $\varepsilon>0$, there exists a constant $N=O(1/\varepsilon^6)$ such that
$|\textrm{cr}(\odot(T^n))/n-c(T)|\leq \varepsilon$ for every integer $n\geq N$.
Let $h:\mathbb{N} \rightarrow \mathbb{R}$ be a function.
In this paper, we gave a necessary and sufficient condition for $\textrm{cr}(\odot(T^n))=h(n)$.
This condition is a corollary of the main result of this paper.
\begin{thm}\label{theo}
Assume that $x_1,...,x_n$ is a list of real numbers and $s :=\sum_{i=1}^{n}x_i$.
Let $h:\mathbb{N} \rightarrow \mathbb{R}$ be a function, and let $c=h(n)/n$. Then
\begin{description}
  \item[(\romannumeral1)]$s<h(n)$ if and only if there exists an integer $k$ such that $\sum_{i=1}^{j}x_{k+i-1}<cj$ for all $j, 1\leq j\leq n$, and $1\leq k\leq n$,
where subscripts are taken modulo $n$.
  \item[(\romannumeral2)]$s>h(n)$ if and only if there exists an integer $k$ such that  $\sum_{i=1}^{j}x_{k+i-1}>cj$ for all $j, 1\leq j\leq n$, and $1\leq k\leq n$,
where subscripts are taken modulo $n$.
\end{description}
\end{thm}
We generalized periodic graphs to generalized periodic graphs (see the definition in Section 2),
which contains periodic graphs and other graphs, such as the complete bipartite graph $K_{m,n}$ and the complete graph $K_n$.
For each $\lambda\in\Lambda$ and a generalized periodic graph $G$, Theorem \ref{theo} implies a necessary and sufficient condition for $\lambda(G)=h(n)$.
Consider, the crossing number, for example.
For a surface $\Sigma$ and a generalized periodic graph $G$,
Theorem \ref{theo} implies a necessary and sufficient condition for $\textrm{cr}_{\Sigma}(G)=h(n)$ (respectively, max-$\textrm{cr}_{\Sigma}(G)=h(n)$).
As an application, we determined the crossing number of the circulant $C(4n;\{1,4\})$.
Theorem \ref{theo} also gives an efficient way to determine whether a generalized periodic graph $G$ can be embedded
in $\Sigma$ if $\Sigma$ is neither $\mathbb{S}_0$ nor $\mathbb{N}_1$.

Base on Theorem \ref{theo}, progress may be made on conjectures below.
Two well-known conjectures on crossing number are Hill's Conjecture,
\begin{equation}\label{Hc}
\textrm{cr}(K_r)=Z(r)Z(r-2)/4,r\geq 3,
\end{equation}
where $Z(r)=\lfloor r/2\rfloor\lfloor(r-1)/2\rfloor$,
and Zarankiewicz's Conjecture,
\begin{equation}\label{Zc}
\textrm{cr}(K_{m,n})=Z(m)Z(n),m,n\geq 1.
\end{equation}
\eqref{Hc} is true \cite{Pan} for $r\leq 12$.
By computer-assisted proofs, McQuillan et al. \cite{McQ} and \'{A}brego et al. \cite{Abr} showed $\textrm{cr}(K_{13})\in\{223,225\}$.
Balogh et al. \cite{Bal} obtained by algorithm that for every sufficiently large $r$, $\textrm{cr}(K_r)>0.985Z(r)Z(r-2)/4$.
Kleitman \cite{Kle} proved that \eqref{Zc} is true for min$(m,n)\leq 6$.
Woodall's algorithm \cite{Woo} confirmed that \eqref{Zc} is true for $K_{7,7}$ and $K_{7,9}$.
Christian et al. \cite{Chr} proved that
for every $m$ there is an $N(m)$ so that if \eqref{Zc} is true for $K_{m,n}$ for all $n\leq N(m)$,
then it is true for all $n$.
Harary et al. \cite{Har} conjectured the crossing number of the Cartesian product $C_m\Box C_n$ is $(m-2)n$, $3\leq m\leq n$, called HKS-Conjecture.
Richter et al. \cite{Ada} verified that this is true for $m\leq 7$.
Glebsky and Salazar \cite{Gle} confirmed that it is true for $m\geq 3$ and $n\geq(m+1)(m+2)/2$.
The crossing number of $C_8\Box C_8$ is unsettled \cite{Chi}.
For progress on the crossing number, readers can refer to \cite{Sch18,Sch22}.
To compute $\textrm{cr}(K_{13})$, $\textrm{cr}(K_{7,11})$ and $\textrm{cr}(C_8\Box C_8)$,
potential algorithms may be devised by Theorem \ref{theo}.

While the focus is on $V$, we defined a partition-transitive graph (see the definition in Section 2).
Let $g:V \rightarrow \mathbb{R}$ be a mapping, and let the \emph{weight} of $g$ on $G$ be $\Sigma_{v\in V}g(v)$.
Let $\Xi$ be the set of parameters that can be expressed as the minimum or maximum of the weight of a specified $g$ on $G$.
For example, the covering number, the stability number and various types of domination number and upper domination number of a graph are members of $\Xi$.
For each $\xi\in\Xi$ and a partition-transitive graph $G$, Theorem \ref{theo} implies a necessary and sufficient condition for $\xi(G)=h(n)$.
Consider, for example, the domination number and the upper domination number.
Let $u$ be a vertex of $G$ and let $N(u)$ be the neighborhood of $u$ in $G$.
A subset $D\subseteq V$ is called a \emph{dominating set} if $N(u)\cap D\not=\emptyset$
for each $u\in V\setminus D$.
The \emph{domination number} $\gamma(G)$ is the minimum cardinality of a dominating set of $G$.
A dominating set $D$ is \emph{minimal} if,
for every vertex $v\in D$, the set $D\setminus\{v\}$ is not a dominating set of $G$.
The \emph{upper domination number} $\Gamma(G)$ is the maximum cardinality of a minimal dominating set of $G$.

Hundreds of types of domination number and upper domination number have been defined.
Many upper bounds on different types of domination number of various classes of graphs \cite{Cho,Spacapan} have been found.
There are many infinite graphs that exact values of various types of domination number \cite{Gon} are known.
Comparatively, exact values of various types of upper domination number of infinite graphs is few,
which are the upper total domination number $\Gamma_t(P_n)$ and the upper paired domination number $\Gamma_p(P_n)$ \cite{Dorbec} of the path $P_n$,
$\Gamma(C_n)$ and $\Gamma_t(C_n)$ \cite{Cyman},
and $\Gamma_p(C_n)$ \cite{Ulatowski}, as far as we known.
For detailed surveys on domination number, we refer the reader to \cite{Haynes20,Haynes21,Haynes23,Henning13}.
As applications, we determined $\gamma_p(C_5\Box C_n)$ and $\Gamma_t(C_4\Box C_n)$.

\section{Generalized periodic graphs and partition-transitive graphs}
Loops and parallel edges are allowed in the graphs of this section.
A periodic graph was introduced by Pinontoan and Richter \cite{Pin04}.
In this paper, we use the definition proposed by Dvo\v{r}\'{a}k and Mohar \cite{RefDvo}.
A \emph{tile} is a triple $Q=(G,L,R),$ where $G$ is a graph, $L=(l_1,...,l_k)$ and
$R=(r_1,...,r_k)$ are sequences of vertices in $G$ of the same length $k$.
Each vertex in $L$ or $R$ can repeatedly appear.
The length of each sequence is called the \emph{width} of the tile.
Suppose that $G_1$ and $G_2$ are two vertex disjoint graphs
and $Q_1=(G_1,L_1,R_1)$ and $Q_2=(G_2,L_2,R_2)$ are tiles of the same width $k$.
Let $Q_1Q_2$ denote the tile $(H,L_1,R_2)$, where $H$ is the graph obtained from the vertex disjoint union of $G_1$ and $G_2$ by adding
an edge between the $j$-th vertices of $R_1$ and $L_2$, $1\leq j\leq k$.
Assume that $Q_1^1=Q_1$ and $Q_1^t=Q_1^{t-1}Q_1$ for integer $t>1$.
Let $Q_1^t=(G_t,L_1,R_t),t\geq 1$.
A \emph{periodic graph} $\odot(Q_1^t),t\geq 2$, is a graph obtained from $G_t$ by adding edges between the $j$-th vertices of $L_1$ and $R_t$ for $j=1,...,k$.
The edges of $\odot(Q_1^t)$ that belong to the copies of $G_1$ are \emph{internal}, and the edges between the copies are \emph{external}.

A \emph{decomposition} of $G$ is a list of edge-disjoint subgraphs of $G$ such that each edge of $G$ appears in one of the subgraphs in the list.
A decomposition $\{H_1,...,H_t\}$ of $G$ is \emph{transitive}
if for each pair of integers $(i,j), i\leq j$,
$H_{i,i+k}$ and $H_{j,j+k}$ are isomorphic for all $k$, $1\leq k\leq t$,
where $H_{i,j}$ is the union of $H_i,...,H_j$ and subscripts are taken modulo $t$.

In a periodic graph $\odot(Q_1^t)$ $(t\geq 3)$, let $S_i$ be the subgraph of $\odot(Q_1^t)$ induced by the set of all external edges between
the $i$-th and the $(i+1)$-th copy of $G_1$ for $1\leq i\leq t-1$,
and let $S_t$ be the subgraph of $\odot(Q_1^t)$ induced by the set of all external edges between the first and the $t$-th copy of $G_1$.
Let $G_i^+$ be the union of $S_i$ and the $i$-th copy of $G_1$ in $\odot(Q_1^t)$ for $1\leq i\leq t$.
It is easy to check that $\{G_1^+,...,G_t^+\}$ is a transitive decomposition of $\odot(Q_1^t)$.

Note that if $1\leq i<j\leq t$ and $V(G_i^+)\cap V(G_j^+)\not=\emptyset$, then either $j-i=1$ or $j-i=t-1$.
We define a similar graph without this restriction.
A graph $G$ is \emph{generalized periodic} if there exists a transitive decomposition of $G$.
Then each periodic graph $\odot(Q_1^t)$ is a generalized periodic graph.
A generalized periodic graph may be not a periodic graph.
Let $U$ and $W$ be the independent sets of $K_{m,n}$, where $U=\{u_i:1\leq i\leq m\}$ and $W=\{v_i:1\leq i\leq n\}$.
Let $V(H_i)=\{u_{i}\}\cup W$ and $E(H_i)=\{u_iv_{j}:1\leq j\leq n\}, 1\leq i\leq m$.
Then $\{H_1,H_2,...,H_{m}\}$ is a transitive decomposition of $K_{m,n}$.
However, $K_{2,3}$ is not a periodic graph, since $|V(K_{2,3})|$ is a prime.

Periodic graphs focus on edges of a graph.
We now turn to vertices of a graph.
A \emph{partition} of $G$ is a list $\{V_1,...,V_t\}$ of subsets of $V(G)$ such that each vertex of $G$ appears in exactly one subset in the list.
A partition $\{V_1,...,V_t\}$ of $G$ is \emph{transitive}
if for each pair of integers $(i,j), i\leq j$,
$V_{i,i+k}$ and $V_{j,j+k}$ are isomorphic for all $k, 1\leq k\leq t$, where $V_{i,j}$ is the union of $V_i,...,V_j$ and subscripts are taken modulo $t$.
A graph $G$ is \emph{partition-transitive} if there exists a transitive partition of $G$.

In a periodic graph $\odot(Q_1^t)$, let $V_i$ be the vertex set of the $i$-th copy of $G_1$, $1\leq i\leq t$.
It is easy to check that $\{V_1,...,V_t\}$ is a transitive partition of $\odot(Q_1^t)$.
Then each periodic graph $\odot(Q_1^t)$ is a partition-transitive graph.
Note that a generalized periodic graph may be not a partition-transitive graph.
For example, $K_{2,3}$ is not a partition-transitive graph.
\section{Bounds on the sum of a list of real numbers and their applications}
\begin{lemma}\label{equalbound}
Suppose that $x_1,...,x_n$ is a list of real numbers and $s :=\sum_{i=1}^{n}x_i$.
Let $h:\mathbb{N} \rightarrow \mathbb{R}$ be a function, and let $c=h(n)/n$. Then
\begin{description}
  \item[(\romannumeral1)]$s\geq h(n)$ if and only if for all $i, 1\leq i\leq n$,
there exists an integer $k_i$ such that $1\leq k_i\leq n$ and $\sum_{j=1}^{k_i}x_{i+j-1}\geq ck_i$,
where subscripts are taken modulo $n$.
  \item[(\romannumeral2)]$s\leq h(n)$ if and only if for all $i, 1\leq i\leq n$,
there exists an integer $k_i$ such that $1\leq k_i\leq n$ and $\sum_{j=1}^{k_i}x_{i+j-1}\leq ck_i$,
where subscripts are taken modulo $n$.
\end{description}
\end{lemma}
\begin{proof}
\textbf{(i)} Suppose that $s\geq h(n)$.
For all $1\leq i\leq n$,
since $\sum_{j=1}^{n}x_{i+j-1}=s$,
$s\geq h(n)$ and $h(n)=cn$, we have
$\sum_{j=1}^{n}x_{i+j-1}\geq cn$, as required.

Conversely, suppose that for all $i, 1\leq i\leq n$, there exists an integer $k_i$ such that $1\leq k_i\leq n$ and $\sum_{j=1}^{k_i}x_{i+j-1}\geq ck_i$.
Let $g^c_i$ be the smallest $k_i$ such that $\sum_{j=1}^{k_i}x_{i+j-1}\geq ck_i$, $1\leq i\leq n$.
Without loss of generality,
assume that max$\{g^c_1,...,g^c_n\}=g^c_1$.
We construct a list of real numbers $y_1,...,y_k$ by the following steps:

\emph{1: set $y_1:=\sum_{i=1}^{g^c_1}x_i$, $k:=1$, $e:=g^c_1$ and $b:=e+1$}

\emph{2: \textbf{while} $e<n$ \textbf{do}}

\emph{3: \quad replace $k$ by $k+1$}

\emph{4: \quad set $e:=b+g^c_b-1$, $y_k:=\sum_{i=b}^ex_i$ and $b:=e+1$}

\emph{5: \textbf{end while}}

\emph{6: return $(y_1,...,y_k)$}

If $g^c_1=n$, then the list of real numbers is $y_1$.
Since $y_1=\sum_{i=1}^nx_i$, we have $\sum_{i=1}^{n}x_i\geq cn$ by the definition of $g^c_1$.
Then $s\geq h(n)$.
Suppose that $g^c_1<n$.
Without loss of generality,
assume that the list of real numbers is $y_1,...,y_k$ with
$y_k=\sum_{i=a}^{a+g^c_a-1}x_i$, where $2\leq a\leq n$.
By the construction of $y_1,...,y_k$, we have $a+g^c_a-1\geq n$.
If $a+g^c_a-1=n$, then $\Sigma_{i=1}^ky_i=\Sigma_{i=1}^nx_i$.
By the construction of $y_1,...,y_k$, we have $\Sigma_{i=1}^ky_i\geq cn$.
Hence $\Sigma_{i=1}^nx_i\geq cn$. Then $s\geq h(n)$.

Otherwise, $a+g^c_a-1>n$.
By the definition of $g^c_a$, we have $\sum_{i=a}^{n}x_i<c(n-a+1)$ and $\sum_{i=a}^{a+g^c_a-1}x_i\geq cg^c_a$.
Therefore, $\sum_{i=n+1}^{a+g^c_a-1}x_i>c(g^c_a-(n-a+1))$.
Since $\sum_{i=n+1}^{a+g^c_a-1}x_i=\sum_{i=1}^{a+g^c_a-1-n}x_i$,
we have $\sum_{i=1}^{a+g^c_a-1-n}x_i>c(g^c_a-(n-a+1))$.
By the definition of $g^c_1$, we have $g^c_a-(n-a+1)\geq g^c_1$.
Hence $g^c_a-g^c_1\geq n-a+1$.
By $a\leq n$, we have $g^c_a-g^c_1\geq 1$.
By max$\{g^c_1,...,g^c_n\}=g^c_1$, we have $g^c_a-g^c_1\leq 0$, a contradiction.

\textbf{(ii)} can be proved similarly.
\end{proof}
Theorem \ref{theo} is the contrapositive of Lemma \ref{equalbound}.
An immediate consequence of Theorem \ref{theo} is the following result.
\begin{cor}\label{coro}
Assume that $x_1,...,x_n$ is a list of real numbers and $0<\varepsilon<1$.
Let $h:\mathbb{N} \rightarrow \mathbb{R}$ be a function.  Then
$\sum_{i=1}^{n}x_i=h(n)$ if and only if there exist integers $k_1$ and  $k_2$ such that $1\leq k_1\leq n$ and $\sum_{i=1}^{j}x_{k_1+i-1}<j(h(n)+\varepsilon)/n$,
and $1\leq k_2\leq n$ and $\sum_{i=1}^{j}x_{k_2+i-1}>j(h(n)-\varepsilon)/n$ for all $j, 1\leq j\leq n$,
where subscripts are taken modulo $n$.
\end{cor}
For each $\lambda\in\Lambda$ and a generalized periodic graph $G$, Theorem \ref{theo} implies a necessary and sufficient condition for $\lambda(G)=h(n)$.
Consider, the domination number and the upper domination number of a graph, for example.
\begin{cor}\label{dom}
Let $h:\mathbb{N} \rightarrow \mathbb{N}$ be a function.
Assume that $G$ is a partition-transitive graph with transitive partition $\{V_1,...,V_n\}$ and $0<\varepsilon<1$.
Let $D$ be a subset of $V(G)$ and let $f_D(V_i)=|V_i\cap D|, 1\leq i\leq n$. Then
$\gamma(G)=h(n)$ if and only if both the following conditions hold.
\begin{description}
  \item[(\romannumeral1)]There exists a dominating set $D$ of $G$ such that for all $j, 1\leq j\leq n$, $\sum_{i=1}^{j}f_D(V_i)<j(h(n)+\varepsilon)/n$.
  \item[(\romannumeral2)]There exists no dominating set $D$ of $G$ such that for all $j, 1\leq j\leq n$, $\sum_{i=1}^{j}f_D(V_i)<j(h(n)-\varepsilon)/n$.
\end{description}
\end{cor}
\begin{proof}
Suppose that $\gamma(G)=h(n)$.
Then there exists a dominating set $D$ of $G$ with $|D|=h(n)$.
Then $f_D(V_1),...,f_D(V_n)$ is a list of $n$ integers and $\sum_{i=1}^{n}f_D(V_i)=h(n)$.
Since $h(n)<h(n)+\varepsilon$, there exists an integer $k$ such that $1\leq k\leq n$ and $\sum_{i=1}^{j}f_D(V_{k+i-1})<j(h(n)+\varepsilon)/n$
for all $j, 1\leq j\leq n$ by Theorem \ref{theo},
where subscripts are taken modulo $n$.
If $k\not=1$, then it is easy to obtain a dominating set $D'$ of $G$ such that $\sum_{i=1}^{j}f_{D'}(V_i)<j(h(n)+\varepsilon)/n$ for all $j, 1\leq j\leq n$,
since $\{V_1,...,V_n\}$ is a transitive partition.
Hence \textbf{(i)} holds.

We prove \textbf{(ii)} by contradiction.
Suppose that there exists a dominating set $D$ of $G$ such that for all $j, 1\leq j\leq n$, $\sum_{i=1}^{j}f_D(V_i)<j(h(n)-\varepsilon)/n$.
Then by Theorem \ref{theo}, $\sum_{i=1}^{n}f_D(V_i)<h(n)-\varepsilon$.
Since $\sum_{i=1}^{n}f_D(V_i)=|D|$, we have $|D|<h(n)-\varepsilon$.
Since $|D|$ is an integer and $\varepsilon>0$, $|D|\leq h(n)-1$.
By the definition of domination number, $|D|\geq h(n)$, a contradiction.

Conversely, suppose that both \textbf{(i)} and \textbf{(ii)} hold.
By \textbf{(i)}, there exists a dominating set $D$ with $|D|<h(n)+\varepsilon$.
Since $|D|$ is an integer and $\varepsilon>0$, $|D|\leq h(n)$.
By the definition of domination number, $\gamma(G)\leq h(n)$.
It suffices to prove $\gamma(G)\geq h(n)$.
By contradiction. Suppose that $\gamma(G)<h(n)$.

Let $D$ be a dominating set of $G$ with $|D|=\gamma(G)$.
Since both $|D|$ and $h(n)$ are integers, $|D|\leq h(n)-1$.
Then $|D|<h(n)-\varepsilon$ by $\varepsilon<1$.
By Theorem \ref{theo}, there exists an integer $k$ such that $1\leq k\leq n$ and $\sum_{i=1}^{j}f_D(V_{k+i-1})<j(h(n)-\varepsilon)/n$
for all $j, 1\leq j\leq n$,
where subscripts are taken modulo $n$.
If $k\not=1$, then it is easy to obtain a dominating set $D'$ of $G$ such that $\sum_{i=1}^{j}f_{D'}(V_i)<j(h(n)-\varepsilon)/n$ for all $j, 1\leq j\leq n$,
since $\{V_1,...,V_n\}$ is a transitive partition.
This contradicts \textbf{(ii)}.
\end{proof}
By a similar argument, Theorem \ref{theo} also implies the following result.
\begin{cor}\label{updom}
Let $h:\mathbb{N} \rightarrow \mathbb{N}$ be a function.
Assume that $G$ is a partition-transitive graph with transitive partition $\{V_1,...,V_n\}$ and $0<\varepsilon<1$.
Let $D$ be a subset of $V(G)$ and let $f_D(V_i)=|V_i\cap D|, 1\leq i\leq n$. Then
$\Gamma(G)=h(n)$ if and only if both the following conditions hold.
\begin{description}
  \item[(\romannumeral1)]There exists a minimal dominating set $D$ of $G$ such that for all $j, 1\leq j\leq n$, $\sum_{i=1}^{j}f_D(V_i)>j(h(n)-\varepsilon)/n$.
  \item[(\romannumeral2)]There exists no minimal dominating set $D$ of $G$ such that for all $j, 1\leq j\leq n$, $\sum_{i=1}^{j}f_D(V_i)>j(h(n)+\varepsilon)/n$.
\end{description}
\end{cor}
Let $D$ be a dominating set of $G$.
A vertex $v$ in $D$ dominates itself and all vertices in $G$ adjacent to it.
Let $u$ be a vertex of $G$, and let $N[u]$ be the closed neighborhood of $u$.
Suppose that $rd_D(u)=|N[u]\cap D|-1$, which counts the number of times a vertex $u$ being re-dominated by $D$.
Let $U\subseteq V(G)$, let $rd_D(U)=\Sigma_{u\in U}rd_D(u)$, and let $rd_D(G)=rd_D(V(G))$.
For the generalized Petersen graph $P(n,2)$,
Fu et al. \cite{Fu09} found that $rd_D(P(n,2))=\Sigma_{u\in V(P(n,2))}(|N[u]\cap D|-1)=4|D|-2n$.
This result can be generalized to the following.
\begin{lemma}\label{rd}
Let $G$ be a $k$-regular graph.
If $D$ is a dominating set of $G$, then $rd_D(G)=(k+1)|D|-|V(G)|$.
\end{lemma}
Lemma  \ref{rd} and Corollary \ref{dom} together imply the following result.
\begin{cor}\label{rddom}
Let $h:\mathbb{N} \rightarrow \mathbb{N}$ be a function.
Assume that $G$ is a partition-transitive graph with transitive partition $\{V_1,...,V_n\}$ and $0<\varepsilon<1$.
Moreover, $G$ is $k$-regular and $t=(k+1)h(n)-|V(G)|$.
Let $D$ be a subset of $V(G)$ and let $rd_D(V_i)=\Sigma_{u\in V_i}rd_D(u), 1\leq i\leq n$,
where $rd_D(u)=|N[u]\cap D|-1$. Then
$\gamma(G)=h(n)$ if and only if both the following conditions hold.
\begin{description}
  \item[(\romannumeral1)]There exists a dominating set $D$ of $G$ such that for all $j, 1\leq j\leq n$, $\sum_{i=1}^{j}rd_D(V_i)<j(t+\varepsilon)/n$.
  \item[(\romannumeral2)]There exists no dominating set $D$ of $G$ such that for all $j, 1\leq j\leq n$, $\sum_{i=1}^{j}rd_D(V_i)<j(t-\varepsilon)/n$.
\end{description}
\end{cor}
For each $\xi\in\Xi$ and a partition-transitive graph $G$, Theorem \ref{theo} implies a necessary and sufficient condition for $\xi(G)=h(n)$.
Consider, the crossing number and the maximum crossing number of a graph in a given surface $\Sigma$, for example.
Let $A$ and $B$ be subgraphs of $G$.
In a drawing $D$ of $G$ in $\Sigma$, the number of crossings crossed by one edge in $A$
and the other edge in $B$ is denoted by $\textrm{cr}_{\Sigma,D}(A,B)$.
Especially, $\textrm{cr}_{\Sigma,D}(A,A)$ is denoted by $\textrm{cr}_{\Sigma,D}(A)$ for short.
The following two equalities are trivial.
\begin{eqnarray}
\textrm{cr}_{\Sigma,D}(A\cup B)&=&\textrm{cr}_{\Sigma,D}(A)+\textrm{cr}_{\Sigma,D}(B)+\textrm{cr}_{\Sigma,D}(A,B),\label{I1}\\
\textrm{cr}_{\Sigma,D}(A,B\cup C)&=&\textrm{cr}_{\Sigma,D}(A,B)+\textrm{cr}_{\Sigma,D}(A,C),\label{I2}
\end{eqnarray}
where $A$, $B$ and $C$ are pairwise edge-disjoint subgraphs of $G$.
The number of crossings of $D$ in $\Sigma$, denoted by $\textrm{cr}_{\Sigma}(D)$, is $\textrm{cr}_{\Sigma,D}(G)$.
Let $H$ be a subgraph of $G$ and let $f_{\Sigma,D}(H): H \rightarrow \mathbb{R}$ be a mapping:
\begin{equation}\label{fun}
f_{\Sigma,D}(H)=\textrm{cr}_{\Sigma,D}(H)+\textrm{cr}_{\Sigma,D}(H,G\setminus E(H))/2.
\end{equation}
By \eqref{I1}, \eqref{I2} and \eqref{fun}, we have:
\begin{lemma}\label{drawing}
Let $\{H_1,...,H_n\}$ be a composition of a graph $G$ and let $D$ be a drawing of $G$ in a surface $\Sigma$.
Then $\textrm{cr}_{\Sigma}(D)=\sum_{i=1}^{n}f_{\Sigma,D}(H_i)$.
\end{lemma}
Lemma \ref{drawing} and Theorem \ref{theo} together imply the following two results.
\begin{cor}\label{cross}
Let $h:\mathbb{N} \rightarrow \mathbb{N}$ be a function.
Assume that $G$ is a generalized periodic graph with transitive decomposition $\{H_1,...,H_n\}$ and $0<\varepsilon<1$.
Let $\Sigma$ be a given surface, let $D$ be a drawing of $G$ in $\Sigma$,
and let $f_{\Sigma,D}(H_i)=\textrm{cr}_{\Sigma,D}(H_i)+\textrm{cr}_{\Sigma,D}(H_i,G\setminus E(H_i))/2, 1\leq i\leq n$. Then
$\textrm{cr}_{\Sigma}(G)=h(n)$ if and only if both the following conditions hold.
\begin{description}
  \item[(\romannumeral1)]There exists a drawing $D$ of $G$ in $\Sigma$ such that $\sum_{i=1}^{j}f_{\Sigma,D}(H_i)<j(h(n)+\varepsilon)/n$ for all $j, 1\leq j\leq n$.
  \item[(\romannumeral2)]There exists no drawing $D$ of $G$ in $\Sigma$ such that $\sum_{i=1}^{j}f_{\Sigma,D}(H_i)<j(h(n)-\varepsilon)/n$ for all $j, 1\leq j\leq n$.
\end{description}
\end{cor}
\begin{cor}\label{maxcross}
Let $h:\mathbb{N} \rightarrow \mathbb{N}$ be a function.
Assume that $G$ is a generalized periodic graph with transitive decomposition $\{H_1,...,H_n\}$ and $0<\varepsilon<1$.
Let $\Sigma$ be a given surface, let $D$ be a drawing of $G$ in $\Sigma$,
and let $f_{\Sigma,D}(H_i)=\textrm{cr}_{\Sigma,D}(H_i)+\textrm{cr}_{\Sigma,D}(H_i,G\setminus E(H_i))/2, 1\leq i\leq n$. Then
max-$\textrm{cr}_{\Sigma}(G)=h(n)$ if and only if both the following conditions hold.
\begin{description}
  \item[(\romannumeral1)]There exists a drawing $D$ of $G$ in $\Sigma$ such that $\sum_{i=1}^{j}f_{\Sigma,D}(H_i)>j(h(n)-\varepsilon)/n$ for all $j, 1\leq j\leq n$.
  \item[(\romannumeral2)]There exists no drawing $D$ of $G$ in $\Sigma$ such that $\sum_{i=1}^{j}f_{\Sigma,D}(H_i)>j(h(n)+\varepsilon)/n$ for all $j, 1\leq j\leq n$.
\end{description}
\end{cor}
\begin{remark}\label{re0}
In Corollary \ref{cross}, if $h(n)=0$, then \rm{\textbf{(i)}} gives an efficient way to determine whether $G$ can be embedded
in $\Sigma$ if $\Sigma$ is neither $\mathbb{S}_0$ nor $\mathbb{N}_1$.
\end{remark}
\begin{remark}\label{re1}
Potential algorithm may be devised to compute $\textrm{cr}(K_{13})$ by Corollary \ref{cross}.
Let $V(K_{13})=\{v_i:1\leq i\leq 13\}$.
Let $V(H_i)=\{v_{i+j}:0\leq j \leq 6\}$ and $E(H_i)=\{v_iv_{i+j}:1\leq j\leq 6\}, 1\leq i\leq 13$, where subscripts are taken modulo $13$.
It is easy to check that $\{H_1,...,H_{13}\}$ is a transitive decomposition of $K_{13}$.
Since $\textrm{cr}(K_{13})\in\{223,225\}$, it suffices to check whether there exists a drawing $D$ of $K_{13}$
with $\textrm{cr}(K_{13})=223$.
By Corollary \ref{cross}, it suffices to check whether there exists a drawing $D$ of $K_{13}$ in the plane such that
$\sum_{i=1}^{j}f_D(H_i)<j(223+\varepsilon)/13$ for all $j, 1\leq j\leq 13$.
Similarly, potential algorithms may be devised to compute $\textrm{cr}(K_{7,11})$ and $\textrm{cr}(C_8\Box C_8)$.
\end{remark}
Let $\textrm{cr}_D$ denote $\textrm{cr}_{\mathbb{S}_0,D}$ for short.
By taking $\Sigma :=\mathbb{S}_0$, the following corollary is immediate from Corollary \ref{cross}.
\begin{cor}\label{periodic}
Let $\{G_1^+,...,G_n^+\}$ be the transitive decomposition of a periodic graph $\odot(Q_1^n)$ as defined in Section 2,
and let $h:\mathbb{N} \rightarrow \mathbb{N}$ be a function.
Assume that $D$ is a drawing of $\odot(Q_1^n)$ in the plane and $0<\varepsilon<1$.
Suppose that $f_D(G_i^+)=\textrm{cr}_D(G_i^+)+\textrm{cr}_D(G_i^+,Q_1^n\setminus E(G_i^+))/2, 1\leq i\leq n$. Then
$\textrm{cr}(Q_1^n)=h(n)$ if and only if both the following conditions hold.
\begin{description}
  \item[(\romannumeral1)]There exists a drawing $D$ of $\odot(Q_1^n)$ in the plane such that for all $j, 1\leq j\leq n$, $\sum_{i=1}^{j}f_D(G_i^+)<j(h(n)+\varepsilon)/n$.
  \item[(\romannumeral2)]There exists no drawing $D$ of $\odot(Q_1^n)$ in the plane such that for all $j, 1\leq j\leq n$, $\sum_{i=1}^{j}f_D(G_i^+)<j(h(n)-\varepsilon)/n$.
\end{description}
\end{cor}

\section{The paired domination number of $C_5\Box C_n$}
Assume that both $m$ and $n$ are integers such that $m\geq 3$ and $n\geq 3$.
Let $V(C_m\Box C_n)=\{u_{i,j}:0\leq i\leq m-1, 0\leq j\leq n-1\}$ and
$E(C_m\Box C_n)=\{u_{i,j}u_{i,j+1},u_{i,j}u_{i+1,j}:0\leq i\leq m-1, 0\leq j\leq n-1\}$, where $i+1$ and $j+1$ are taken modulo $m$ and $n$,
respectively.
Assume that $V_j=\{u_{i,j}:0\leq i\leq m-1\}, 0\leq j\leq n-1$.
It is easy to check that $\{V_0,...,V_{n-1}\}$ is a transitive partition of $C_m\Box C_n$.
Let $S$ be a subset of $V(C_m\Box C_n)$.
Assume that $V_{i,j},i\leq j,$ is the union of $V_i,...,V_j$ with subscripts taken modulo $n$,
and $f_S(V_{i,j})=|V_{i,j}\cap S|$.

A \emph{matching} in a graph is a set of pairwise nonadjacent edges.
A matching $M$ of graph $G$ is \emph{perfect} if each vertex of $G$ incident with an edge of $M$.
A dominating set $D$ of graph $G$ is \emph{paired} if the induced subgraph $G[D]$ of $G$ contains a perfect matching \cite{Haynes95,Haynes98}.
The \emph{paired domination number} of graph $G$, denoted by $\gamma_p(G)$, is the minimum cardinality of a paired dominating set of $G$.
For $\gamma_p(C_m\Box C_n)$, the exact values \cite{Hu14} were obtained for $m=3,4$,
and an upper bound \cite{Hu16} was given for $m=5$.
Hu and Xu proved the following result.
\begin{lemma}[\cite{Hu14}]\label{pair34}
$\gamma_p(C_5\Box C_n)=\lceil 4n/3\rceil, n=3,4$.
\end{lemma}
Let
\[h(n)=\left\{
\begin{aligned}
\lceil 4n/3\rceil+1,\;&\text{\emph{if}}\; n\equiv 2\,(\rm{mod}\; 3);\\
\lceil 4n/3\rceil,\quad\ \; \;&\text{\emph{otherwise}}.
\end{aligned}
\right.\]
Hu et al. gave the following bound.
\begin{lemma}[\cite{Hu16}]\label{uppair}
$\gamma_p(C_5\Box C_n)\leq h(n), n\geq 3$.
\end{lemma}
Haynes and Slater gave the following result.
\begin{lemma}[\cite{Haynes98}]\label{gb}
Suppose that $G$ is a graph with no isolated vertices, $|V(G)|=n$ and maximum degree $\Delta$.
Then $\gamma_p(G)\geq n/\Delta$, and this bound is sharp.
\end{lemma}
The degree of a vertex $v$ in a graph $G$, denoted by $d_G(v)$, is the number of edges of $G$ incident with $v$.
By the definition of a paired dominating set, the following result is trivial.
\begin{lemma}\label{pairset}
Let $S$ be a paired dominating set of a graph $G$ and let $H$ be a component of the induced subgraph $G[S]$ of $G$.
Then:
\begin{description}
  \item[(\romannumeral1)]both $|S|$ and $|V(H)|$ are even integers,
  \item[(\romannumeral2)]for each $u\in V(H)$ with $d_H(u)\geq 3$, $u$ has at most one neighbor $v$ in $H$ such that $d_H(v)=1$.
\end{description}
\end{lemma}
\begin{thm}\label{t1}
$\gamma_p(C_5\Box C_n)=h(n), n\geq 3$.
\end{thm}
\begin{proof}
It suffices to prove $\gamma_p(C_5\Box C_n)\geq h(n)$ by Lemma \ref{uppair}.
We prove it by induction on $n$.
By Lemma \ref{gb}, we have $\gamma_p(C_5\Box C_5)\geq 25/4$.
Then $\gamma_p(C_5\Box C_5)\geq 8$, since $\gamma_p(C_5\Box C_5)$ is even by Lemma \ref{pairset}.
Combining Lemma \ref{pair34}, it is true for $3\leq n\leq 5$.
Assume that $n\geq 5$ in the following of this proof.
Suppose that $\gamma_p(C_5\Box C_n)\geq h(n)$
and $\gamma_p(C_5\Box C_r)\geq h(r)$, $3\leq r\leq n$.
It suffices to prove $\gamma_p(C_5\Box C_{n+1})\geq h(n+1)$.
By contradiction. Suppose that $\gamma_p(C_5\Box C_{n+1})<h(n+1)$.
Let $S$ be a paired dominating set of $C_5\Box C_{n+1}$ with $|S|=\gamma_p(C_5\Box C_{n+1})$.

If $(n+1)\equiv 2$ (mod 3), then $|S|<\lceil 4(n+1)/3\rceil+1$.
Since both $|S|$ and $\lceil 4(n+1)/3\rceil+1$ are even integers, we have $|S|<\lceil 4(n+1)/3\rceil$.
Hence $|S|<\lceil 4(n+1)/3\rceil$ for all $n$ by $|S|<h(n+1)$.
Moreover, since $|S|$ is an integer, $|S|<4(n+1)/3$.
By Lemma \ref{rd}, $rd_S(C_5\Box C_{n+1})=5|S|-5(n+1)$.
Hence $rd_S(C_5\Box C_{n+1})<5(n+1)/3$.
Then $rd_S(V_0),...,rd_S(V_n)$ is a list of $n+1$ integers,
and $\Sigma_{i=0}^nrd_S(V_i)=rd_S(C_5\Box C_{n+1})$.
By Theorem \ref{theo}, without loss of generality, assume that $rd_S(V_{0,i})<5(i+1)/3$ for all $i, 0\leq i\leq n$.
Then $rd_S(V_0)\leq 1$, $rd_S(V_{0,1})\leq 3$, $rd_S(V_{0,2})\leq 4$, $rd_S(V_{0,3})\leq 6$, $rd_S(V_{0,4})\leq 8$ and $rd_S(V_{0,5})\leq 9$.
Let $\overline{S}=V(C_5\Box C_{n+1})\setminus S$.
Since $S$ is a paired dominating set of $C_5\Box C_{n+1}$,
$|S|$ is even by Lemma \ref{pairset}.
Moreover,
if $u\in S$, then $rd_S(u)\geq 1$.
By $rd_S(V_0)\leq 1$, $|V_0\cap S|\leq 1$.
We consider two cases.

Case 1. $|V_0\cap S|=1$.

Without loss of generality, assume that $u_{2,0}\in S$.
Then $rd_S(u_{2,0})\geq 1$.
By $rd_S(V_0)\leq 1$, $\{u_{1,1}, u_{3,1}\}\subseteq \overline{S}$.
By $rd_S(V_{0,1})\leq 3$, $|V_1\cap S|\leq 2$.
If $|V_1\cap S|=0$, then $\{u_{0,2}, u_{1,2}, u_{3,2}, u_{4,2}\}\subseteq S$,
since each of $u_{0,1}$, $u_{1,1}$, $u_{3,1}$ and $u_{4,1}$ is dominated by $S$.
Hence $rd_S(V_{0,2})\geq 5$, a contradiction.
We consider the following subcases:

Subcase 1.1. $|V_1\cap S|=1$.

By symmetry, either $u_{0,1}$ or $u_{2,1}$ is in $S$.
If $u_{0,1}\in S$, then $\{u_{0,2},u_{3,2}\}\subseteq S$,
since both $u_{0,1}$ and $u_{3,1}$ are dominated by $S$.
Hence $\{u_{2,0}, u_{0,1}, u_{0,2}, u_{3,2}\}\subseteq S$ and $rd_S(u_{4,2})\geq 1$.
Then $rd_S(V_{0,2})\geq 5$, a contradiction.

If $u_{2,1}\in S$, then $\{u_{0,2},u_{4,2}\}\subseteq S$,
since both $u_{0,1}$ and $u_{4,1}$ are dominated by $S$.
By $rd_S(V_{0,2})\leq 4$, $\{u_{1,2},u_{2,2},u_{3,2}\}\cup V_3\subseteq \overline{S}$.
Since each of $u_{1,3}$, $u_{2,3}$ and $u_{3,3}$ is dominated by $S$,
$\{u_{1,4}, u_{2,4}, u_{3,4}\}\subseteq S$.
By $rd_S(V_{0,4})\leq 8$ and $rd_S(u_{2,4})\geq 2$, $\{u_{0,4},u_{4,4}\}\cup V_5\subseteq\overline{S}$.
Then $G[S]$ has a component $G[\{u_{1,4},u_{2,4},u_{3,4}\}]$ with 3 vertices, which contradicts Lemma \ref{pairset}.

Subcase 1.2. $|V_1\cap S|=2$.

By symmetry, either $\{u_{0,1},u_{2,1}\}\subseteq S$ or $\{u_{0,1},u_{4,1}\}\subseteq S$.
If $\{u_{0,1},u_{2,1}\}\subseteq S$, then since $rd_S(u_{1,1})\geq 1$, we have $rd_S(V_{0,1})\geq 4$, a contradiction.
Suppose that $\{u_{0,1},u_{4,1}\}\subseteq S$.
By $rd_S(V_{0,1})\leq 3$, $\{u_{2,1}\}\cup V_2\subseteq \overline{S}$.
Since each of $u_{1,2}$, $u_{2,2}$ and $u_{3,2}$ is dominated by $S$,
$\{u_{1,3}, u_{2,3}, u_{3,3}\}\subseteq S$.
Then $rd_S(u_{2,3})\geq 2$, and $rd_S(V_{0,3})\geq 7$, a contradiction.

Case 2. $|V_0\cap S|=0$.

By $rd_S(V_{0,1})\leq 3$, $|V_1\cap S|\leq 3$.
If $|V_1\cap S|=0$, then $V_2\subseteq S$, since each vertex in $V_1$ is dominated by $S$.
Hence $rd_S(V_{0,2})\geq 5$, a contradiction.
Suppose that $|V_1\cap S|=1$. Without loss of generality, assume that $u_{2,1}\in S$.
Since each of $u_{0,1}$, $u_{2,1}$ and $u_{4,1}$ is dominated by $S$,
$\{u_{0,2}, u_{2,2}, u_{4,2}\}\subseteq S$.
Then $rd_S(u_{1,2})\geq 1$.
Hence $rd_S(V_{0,2})\geq 5$, a contradiction.
Suppose that $|V_1\cap S|=3$. Then without loss of generality, assume that $\{u_{0,1},u_{4,1}\}\subseteq S$.
By symmetry, either $u_{1,1}\in S$ or $u_{2,1}\in S$.
If $u_{1,1}\in S$, then $rd_S(u_{0,1})\geq 2$;
if $u_{2,1}\in S$, then $rd_S(u_{1,1})\geq 1$.
In both cases, $rd_S(V_{0,1})\geq 4$, a contradiction.

Suppose that $|V_1\cap S|=2$.
Then either $\{u_{0,1},u_{4,1}\}\subseteq S$ or $\{u_{1,1},u_{3,1}\}\subseteq S$ by symmetry.
If $\{u_{1,1},u_{3,1}\}\subseteq S$, then $\{u_{1,2},u_{3,2}\}\subseteq S$,
since both $u_{1,1}$ and $u_{3,1}$ are dominated by $S$.
Then $rd_S(u_{2,1})\geq 1$, and $rd_S(V_{0,2})\geq 5$, a contradiction.

Suppose that $\{u_{0,1},u_{4,1}\}\subseteq S$.
Since $u_{2,1}$ is dominated by $S$, $u_{2,2}\in S$.
By $rd_S(V_{0,2})\leq 4$, $\{u_{0,2},u_{1,2},u_{3,2},u_{4,2}\}\subseteq\overline{S}$.
Since $u_{2,2}$ is dominated by $S$, $u_{2,3}\in S$.
By $rd_S(V_{0,3})\leq 6$, $\{u_{0,3},u_{1,3},u_{3,3},u_{4,3}\}\subseteq\overline{S}$.
Since both $u_{0,3}$ and $u_{4,3}$ are dominated by $S$,
$\{u_{0,4},u_{4,4}\}\subseteq S$.
By $rd_S(V_{0,4})\leq 8$, $\{u_{1,4},u_{2,4},u_{3,4}\}\subseteq\overline{S}$.
By $rd_S(V_{0,5})\leq 9$, $|V_5\cap S|\leq 1$. We consider the following subcases.

Subcase 2.1. $|V_5\cap S|=1$.

If $n=5$, then $|S|$ is odd, a contradiction.
Suppose that $n\geq 6$.
Then $rd_S(V_{0,6})<35/3$.
By symmetry, one of $u_{0,5},u_{1,5}$ and $u_{2,5}$ in in $S$.
If $u_{0,5}\in S$, then $\{u_{2,6},u_{3,6}\}\subseteq S$,
since both $u_{2,5}$ and $u_{3,5}$ are dominated by $S$.
By Lemma \ref{pairset}, $u_{0,6}\in S$.
Hence $rd_S(u_{0,4})\geq 2$ and $rd_S(u_{0,5})\geq 2$.
If $u_{1,5}\in S$, then $\{u_{1,6},u_{3,6}\}\subseteq S$,
since both $u_{1,5}$ and $u_{3,5}$ are dominated by $S$.
Hence $rd_S(u_{1,4})\geq 1$, $rd_S(u_{0,5})\geq 1$ and $rd_S(u_{2,6})\geq 1$.
In both cases, $rd_S(V_{0,6})\geq 12$, a contradiction.

Suppose that $u_{2,5}\in S$.
Let $G_1=C_5\Box C_{n+1}-V_{2,4}$, let $G_2$ be the graph with $V(G_2)=V(G_1)$ and let
$E(G_2)=E(G_1)\cup \{u_{i,1}u_{i,5}: 0\leq i\leq 4\}$.
It is easy to check that $G_2$ is isomorphic to $C_5\Box C_{n-2}$
and $S\setminus V_{2,4}$ is a paired dominating set of $G_2$.
Since $\gamma_p(C_5\Box C_{n-2})\geq h(n-2)$ and $|V_{2,4}\cap S|=4$,
we have $|S\setminus V_{2,4}|\geq h(n-2)$ and $|S\setminus V_{2,4}|=|S|-4$, respectively.
Then $|S|=|S\setminus V_{2,4}|+4\geq h(n-2)+4\geq h(n+1)$,
a contradiction.

Subcase 2.2. $|V_5\cap S|=0$.

If $n=5$, then $u_{2,5}$ is not dominated by $S$, a contradiction.
Suppose that $n\geq 6$.
Then $rd_S(V_{0,6})<35/3$.
Since each of $u_{1,5}$, $u_{2,5}$ and $u_{3,5}$ is dominated by $S$,
$\{u_{1,6},u_{2,6},u_{3,6}\}\subseteq S$.
Since $rd_S(u_{2,6})\geq 2$, $\{u_{0,6},u_{4,6}\}\subseteq\overline{S}$ by $rd_S(V_{0,6})<35/3$.
If $n=6$, then $|S|$ is odd, a contradiction.
Suppose that $n\geq 7$.
Then $rd_S(V_{0,7})<40/3$.
By $rd_S(V_{0,7})<40/3$, $|V_7\cap S|\leq 1$.
Then by Lemma \ref{pairset}, assume that $u_{1,7}\in S$ and $\{u_{0,7},u_{2,7},u_{3,7},u_{4,7}\}\subseteq\overline{S}$ by symmetry.
Then $rd_S(u_{1,6})\geq 2$ and $rd_S(u_{2,7})\geq 1$.

If $n=7$, then $u_{4,7}$ is not dominated by $S$, a contradiction.
Suppose that $n\geq 8$.
Then $rd_S(V_{0,8})<45/3$.
Since $u_{4,7}$ is dominated by $S$, $u_{4,8}\in S$ and $\{u_{0,8},u_{1,8},u_{2,8},u_{3,8}\}\subseteq\overline{S}$ by $rd_S(V_{0,8})<45/3$.
If $n=8$, then $|S|$ is odd, a contradiction.
Suppose that $n\geq 9$.
Then $rd_S(V_{0,9})<50/3$.
Since both $u_{2,8}$ and $u_{4,8}$ are dominated by $S$,
$\{u_{2,9},u_{4,9}\}\subseteq S$.
Then $rd_S(u_{3,9})\geq 1$, and $rd_S(V_{0,9})\geq 17$, a contradiction.
\end{proof}

\section{The upper total domination number of $C_4\Box C_n$}
A \emph{total dominating set}, abbreviated TD-set, of a graph $G$ with no isolated vertex is a subset $S$ of $V(G)$
such that every vertex in $G$ is adjacent to a vertex in $S$.
The \emph{total domination number} \cite{Cockayne} of $G$,
denoted by $\gamma_t(G)$, is the minimum cardinality of a TD-set of $G$.
A TD-set $S$ of $G$ is \emph{minimal} if no proper subset of $S$ is a TD-set of $G$.
The \emph{upper total domination number} $\Gamma_t(G)$ of $G$ is the maximum cardinality of a minimal TD-set of $G$.
As mentioned earlier, only both $\Gamma_t(P_n)$ \cite{Dorbec} and $\Gamma_t(C_n)$ \cite{Cyman} were obtained, as far as we known.
Moreover, exact values of $\gamma_t(C_m\Box C_n)$ \cite{Hu16,Hu14} were determined for $m=3,4,5$.

Let $G=(V,E)$ be a graph.
For a set $S\subseteq V$ and a vertex $v\in S$,
the open $S$-private neighborhood of $v$, denoted by pn$(v,S)$, is the set $N(v)\setminus N(S\setminus\{v\})$.
Equivalently, pn$(v,S)=\{w\in V:N(w)\cap S=\{v\}\}$.
Let the open $S$-external private neighborhood epn$(v,S)$ of $v$ be pn$(v,S)\setminus S$.
The open $S$-internal private neighborhood of $v$ is defined by ipn$(v,S)=$ pn$(v,S)\cap S$.
Note that pn$(v,S)=$ ipn$(v,S)\;\cup$ epn$(v,S)$.
Cockayne et al. proved the following result.
\begin{lemma}[\cite{Cockayne}]\label{n1}
A TD-set $S$ of graph $G$ is a minimal TD-set of $G$
if and only if ipn$(v,S)\ne\emptyset$ or epn$(v,S)\ne\emptyset$ for every vertex $v\in S$.
\end{lemma}
\begin{thm}\label{n4}
$\Gamma_t(C_{4}\Box C_{n})=2n$, $n\geq 3$.
\end{thm}
\begin{proof}
Since $\{u_{1,i},u_{2,i}:0\leq i\leq n-1\}$ is a minimal TD-set of $C_{4}\Box C_{n}$,
we have $\Gamma_t(C_{4}\Box C_{n})\geq 2n$.
It suffices to prove $\Gamma_t(C_{4}\Box C_n)\leq 2n$.
By contradiction. Suppose that $\Gamma_t(C_{4}\Box C_n)>2n$.
Let $S$ be a minimal TD-set of $C_4\Box C_n$ with $|S|=\Gamma_t(C_4\Box C_n)$.
Then $f_S(V_0),...,f_S(V_{n-1})$ is a list of $n$ integers,
and $\Sigma_{i=0}^{n-1}f_S(V_i)=|S|$.
By Theorem \ref{theo}, without loss of generality, assume that $f_S(V_{0,i})>2(i+1), 0\leq i\leq n-1$.
Hence $f_S(V_0)\geq 3$, $f_S(V_{0,1})\geq 5$ and $f_S(V_{0,2})\geq 7$,
and $f_S(V_{0,3})\geq 9$ if $n\geq 4$.
Let $\overline{S}=V(C_4\Box C_n)\setminus S$.
We distinguish two cases:

Case 1. $f_S(V_0)=4$.

Then $V_0\subseteq S$.
By $f_S(V_{0,1})\geq 5$, we have $|V_1\cap S|\geq 1$.
Without loss of generality, assume that $u_{3,1}\in S$.
We consider the following subcases.

Subcase 1.1. $u_{3,2}\in \overline{S}$.

By ipn$(u_{3,1},S)=\emptyset$ and Lemma \ref{n1}, epn$(u_{3,1},S)\neq\emptyset$.
If $n=3$, then epn$(u_{3,1},S)=\emptyset$, a contradiction.
Suppose that $n\geq 4$.
Then by epn$(u_{3,1},S)\neq\emptyset$, $\{u_{0,2},u_{2,2},u_{3,3}\}\subseteq \overline{S}$.

Suppose that $u_{2,1}\in S$.
By ipn$(u_{2,1},S)=\emptyset$ and Lemma \ref{n1}, epn$(u_{2,1},S)\neq\emptyset$.
Hence $\{u_{1,2},u_{2,3}\}\subseteq \overline{S}$.
By $f_S(V_{0,2})\geq 7$, without loss of generality, assume that $u_{1,1}\in S$.
By ipn$(u_{1,1},S)=\emptyset$ and Lemma \ref{n1}, epn$(u_{1,1},S)\neq\emptyset$.
Hence $u_{1,3}\in\overline{S}$.
By $f_S(V_{0,3})\geq 9$, $\{u_{0,1},u_{0,3}\}\subseteq S$.
Then ipn$(u_{0,1},S)=$ epn$(u_{0,1},S)=\emptyset$, a contradiction.
Hence $u_{2,1}\in\overline{S}$.
By symmetry, $u_{0,1}\in\overline{S}$.

By $f_S(V_{0,2})\geq 7$, $\{u_{1,1},u_{1,2}\}\subseteq S$.
By epn$(u_{1,1},S)=\emptyset$ and Lemma \ref{n1}, ipn$(u_{1,1},S)\neq\emptyset$.
Hence $u_{1,3}\in \overline{S}$.
By $f_S(V_{0,3})\geq 9$, $\{u_{0,3},u_{2,3}\}\subseteq S$.
Then ipn$(u_{1,2},S)=$ epn$(u_{1,2},S)=\emptyset$, a contradiction.

Subcase 1.2. $u_{3,2}\in S$.

By epn$(u_{3,1},S)=\emptyset$ and Lemma \ref{n1}, ipn$(u_{3,1},S)\neq\emptyset$.
If $n=3$, then ipn$(u_{3,1},S)=\emptyset$, a contradiction.
Suppose that $n\geq 4$.
Then $\{u_{0,2},u_{2,2},u_{3,3}\}\subseteq \overline{S}$ by ipn$(u_{3,1},S)\neq\emptyset$.
If $u_{2,1}\in S$, then ipn$(u_{2,1},S)=$ epn$(u_{2,1},S)=\emptyset$, a contradiction.
Hence $u_{2,1}\in \overline{S}$.
By symmetry, $u_{0,2}\in \overline{S}$.
By ipn$(u_{3,2},S)=\emptyset$ and Lemma \ref{n1}, epn$(u_{3,2},S)\neq\emptyset$.
If $u_{3,3}\in$ epn$(u_{3,2},S)$, then $\{u_{0,3},u_{2,3}\}\subseteq \overline{S}$.
By $f_S(V_{0,3})\geq 9$, $\{u_{1,1},u_{1,2},u_{1,3}\}\subseteq S$.
Otherwise, assume that $u_{2,2}\in$ epn$(u_{3,2},S)$ by symmetry.
Then $\{u_{1,2},u_{2,3}\}\subseteq \overline{S}$.
By $f_S(V_{0,3})\geq 9$, $\{u_{1,1},u_{0,3},u_{1,3}\}\subseteq S$.
In both cases, ipn$(u_{1,1},S)=$ epn$(u_{1,1},S)=\emptyset$, a contradiction.

Case 2. $f_S(V_0)=3$.

By symmetry, assume that $\{u_{0,0},u_{1,0},u_{2,0}\}\subseteq S$ and $u_{3,0}\in\overline{S}$.
By $f_S(V_{0,1})\geq 5$, $|V_1\cap S|\geq 2$.
We consider the following subcases.

Subcase 2.1. $u_{3,1}\in S$.

Suppose that $u_{2,1}\in S$.
If $u_{3,2}\in S$, then ipn$(u_{2,1},S)=$ epn$(u_{2,1},S)=\emptyset$, a contradiction.
Hence $u_{3,2}\in \overline{S}$.
If $n=3$, then ipn$(u_{2,0},S)=$ epn$(u_{2,0},S)=\emptyset$, a contradiction.
Suppose that $n\geq 4$.
By ipn$(u_{3,1},S)=\emptyset$ and Lemma \ref{n1}, epn$(u_{3,1},S)\neq\emptyset$.
Then $\{u_{0,2},u_{2,2},u_{3,3}\}\subseteq \overline{S}$.

If $u_{0,1}\in S$, then by ipn$(u_{0,1},S)=\emptyset$ and Lemma \ref{n1}, epn$(u_{0,1},S)\neq\emptyset$.
Hence $\{u_{1,2},u_{0,3}\}\subseteq \overline{S}$.
By $f_S(V_{0,3})\geq 9$, $\{u_{1,1},u_{1,3},u_{2,3}\}\subseteq S$.
Then ipn$(u_{1,1},S)=$ epn$(u_{1,1},S)=\emptyset$, a contradiction.
Hence $u_{0,1}\in \overline{S}$.
By $f_S(V_{0,2})\geq 7$, $\{u_{1,1},u_{1,2}\}\subseteq S$.
By epn$(u_{1,1},S)=\emptyset$ and Lemma \ref{n1}, ipn$(u_{1,1},S)\neq\emptyset$.
Hence $u_{1,3}\in \overline{S}$.
By $f_S(V_{0,3})\geq 9$, $\{u_{0,3},u_{2,3}\}\subseteq S$.
Then ipn$(u_{1,2},S)=$ epn$(u_{1,2},S)=\emptyset$, a contradiction.

Therefore, $u_{2,1}\in \overline{S}$. By symmetry, $u_{0,1}\in \overline{S}$.
By $f_S(V_{0,1})\geq 5$, $u_{1,1}\in S$.
Since $u_{3,1}$ is dominated by $S$, $u_{3,2}\in S$.
If $n=3$, then ipn$(u_{1,1},S)=$ epn$(u_{1,1},S)=\emptyset$, a contradiction.
Suppose that $n\geq 4$.
By epn$(u_{3,1},S)=\emptyset$ and Lemma \ref{n1}, ipn$(u_{3,1},S)\neq\emptyset$.
Then $\{u_{0,2},u_{2,2},u_{3,3}\}\subseteq \overline{S}$.
By $f_S(V_{0,2})\geq 7$, $u_{1,2}\in S$.
By epn$(u_{1,1},S)=\emptyset$ and Lemma \ref{n1}, ipn$(u_{1,1},S)\neq\emptyset$.
Then $u_{1,3}\in\overline{S}$.
By $f_S(V_{0,3})\geq 9$, $\{u_{0,3},u_{2,3}\}\subseteq S$.
Then ipn$(u_{1,2},S)=$ epn$(u_{1,2},S)=\emptyset$, a contradiction.

Subcase 2.2. $u_{3,1}\in \overline{S}$.

By $f_S(V_{0,1})\geq 5$, we consider the following cases.

Subcase 2.2.1. $f_S(V_{0,1})=6$.

If $u_{0,2}\in S$, then ipn$(u_{1,1},S)=$ epn$(u_{1,1},S)=\emptyset$, a contradiction.
Hence $u_{0,2}\in\overline{S}$.
By symmetry, $u_{2,2}\in \overline{S}$.
By ipn$(u_{0,1},S)=\emptyset$ and Lemma \ref{n1}, epn$(u_{0,1},S)\neq\emptyset$.
Then $\{u_{1,2},u_{3,2}\}\subseteq \overline{S}$.
Then $f_S(V_{0,2})=6$, a contradiction.

Subcase 2.2.2. $f_S(V_{0,1})=5$.

By symmetry, either $\{u_{0,1},u_{2,1}\}\subseteq S$ or $\{u_{1,1},u_{2,1}\}\subseteq S$.
Suppose that $\{u_{0,1},u_{2,1}\}\subseteq S$.
If $u_{0,2}\in\overline{S}$, then by ipn$(u_{0,1},S)=\emptyset$ and Lemma \ref{n1}, epn$(u_{0,1},$ $S)\neq\emptyset$.
Hence $\{u_{1,2},u_{3,2}\}\subseteq \overline{S}$ and $f_S(V_{0,2})\leq 6$, a contradiction.
Then $u_{0,2}\in S$.
By symmetry, $u_{2,2}\in S$.
If $n=3$, then ipn$(u_{1,0},S)=$ epn$(u_{1,0},S)=\emptyset$, a contradiction.
Suppose that $n\geq 4$.
By epn$(u_{0,1},S)=\emptyset$ and Lemma \ref{n1}, ipn$(u_{0,1},S)\neq\emptyset$.
Hence $\{u_{1,2},u_{3,2},u_{0,3}\}\subseteq\overline{S}$.
By symmetry, $u_{2,3}\in \overline{S}$.
By $f_S(V_{0,3})\geq 9$, $\{u_{1,3},u_{3,3}\}\subseteq S$.
Then ipn$(u_{0,2},S)=$ epn$(u_{0,2},S)=\emptyset$, a contradiction.

Suppose that $\{u_{1,1},u_{2,1}\}\subseteq S$.
If $u_{1,2}\in\overline{S}$, then by ipn$(u_{1,1},S)=\emptyset$ and Lemma \ref{n1}, epn$(u_{1,1},S)\neq\emptyset$.
Hence $\{u_{0,2},u_{2,2}\}\subseteq \overline{S}$ and $f_S(V_{0,2})\leq 6$, a contradiction.
Suppose that $u_{1,2}\in S$.
Then by epn$(u_{1,1},S)=\emptyset$ and Lemma \ref{n1}, ipn$(u_{1,1},S)\neq\emptyset$.
Hence $\{u_{0,2},u_{2,2}\}\subseteq \overline{S}$.
By $f_S(V_{0,2})\geq 7$, $u_{3,2}\in S$.
Then ipn$(u_{2,1},S)=$ epn$(u_{2,1},S)=\emptyset$, a contradiction.
\end{proof}

\section{The crossing number of $C(4k;\{1,4\})$}
\begin{figure}
  \centering
  \subfigure[$C(4k;\{1,4\})$]{
    \label{fig:subfig:a}
    \includegraphics[width=0.4\textwidth]{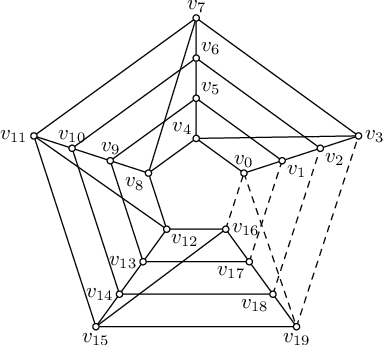}}
    \subfigure[A subdivision of $C(4(k-1);\{1,4\})$]{
    \label{fig:subfig:a}
    \includegraphics[width=0.4\textwidth]{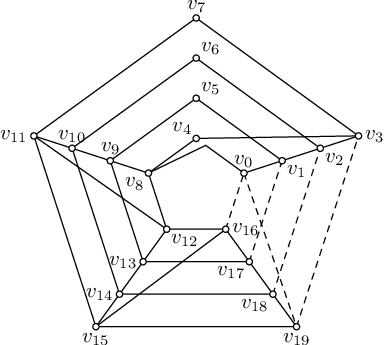}}
  \hspace{0.1in}
  \caption{Good drawings of $C(4k;\{1,4\})$ and a subdivision of $C(4(k-1);\{1,4\})$.}\label{F1}
\end{figure}

In this section, we consider the crossing number of a graph in the plane.
Let $n$ be a positive integer and let $\mathbb{Z}_n=\{0,...,n-1\}$.
A circulant is a Cayley graph for a cyclic group.
To be more precise, if $S$ is any subset of $\mathbb{Z}_n\setminus \{0\}$,
then a circulant $C(n;S)$ is an undirected simple graph of order $n$ with vertex-set $\mathbb{Z}_n$
and edge-set $\{\{i,i+a\}:i\in \mathbb{Z}_n, a\in S\}$.
Yang and Zhao \cite{RefYang01} calculated
$\textrm{cr}(C(n;\{1,k\}))$ for $5\leq n\leq 18$ and $2\leq k\leq \lfloor n/2\rfloor$ by the CCN algorithm \cite{RefLin09}.
Yang et al. \cite{RefYang04} proved $\textrm{cr}(C(n;\{1,3\}))=\lfloor n/3\rfloor+(n$ mod $3$), $n\geq 8$.
Ma et al. \cite{RefMa05} showed $\textrm{cr}(C(2k+2;\{1,k\}))=k+1$, $k\geq 3$.
Lin et al. \cite{RefLin05} proved $\textrm{cr}(C(3k;\{1,k\}))=k$, $k\geq 3$,
from which $\textrm{cr}(C(12;\{1,4\}))=4$.
Ho \cite{RefHo} showed $\textrm{cr}(C(3k+1;\{1,k\}))=k+1$, $k\geq 3$.
Since $C(8;\{1,4\})$ is isomorphic to the M\"{o}bius ladder of order 8 and the crossing number of each M\"{o}bius ladder \cite{RefGH} is 1,
we have $\textrm{cr}(C(8;\{1,4\}))=1$.
Yang and Zhao gave the following result.
\begin{lemma}[\cite{RefYang01}]\label{C16}
$\textrm{cr}(C(16;\{1,4\}))=8$.
\end{lemma}
In this section, our main result is the following theorem.
\begin{thm}\label{T41}
$\textrm{cr}(C(4k;\{1,4\}))=2k$, $k\geq 4$.
\end{thm}
Let $k\geq 4$, let $V(C(4k;\{1,4\}))=\{v_i:i\in\mathbb{Z}_{4k}\}$ and let
$E(C(4k;\{1,4\}))=\{v_iv_{i+a}:i\in\mathbb{Z}_{4k},a\in\{1,4\}\}$, where subscripts are taken modulo $4k$.
By the drawing in Fig. \ref{F1}a, we have the following lemma.
\begin{lemma}\label{lemma41}
$\textrm{cr}(C(4k;\{1,4\}))\leq 2k$, $k\geq 4$.
\end{lemma}
In the following of this section, assume that $k\geq 5$, $i\in\mathbb{Z}_{4k}$ and $G^{4k}=C(4k;\{1,4\})$ if not specified.
Let $V(H_i)=\{v_{i},v_{i+1},v_{i+4}\}$ and $E(H_i)=\{v_{i}v_{i+1},$ $v_{i}v_{i+4}\}$.
It is easy to check that $\{H_0,H_1,...,H_{4k-1}\}$ is a transitive decomposition of $G^{4k}$.
Let
$C(v_i,v_{i+1},v_{i+5},v_{i+4})$ be the cycle $v_iv_{i+1}v_{i+5}v_{i+4}v_i$,
and let $P(v_{i+1},v_{i+2},v_{i+3},v_{i+4})$ be the path $v_{i+1}v_{i+2}v_{i+3}v_{i+4}$.
$C(v_i,v_{i+1},v_{i+2},v_{i+3},v_{i+4})$ and $P(v_i,v_{i+1},v_{i+2},v_{i+3},v_{i+4})$ are defined similarly.

A drawing of a graph in the plane with the minimum number of pairwise crossings of edges is always a good drawing.
A drawing of a graph is \emph{good}, provided that no edges crosses itself,
no adjacent edges cross each other, no two edges cross more than once,
and no three edges cross at a common point.
In a good drawing, a common point of two edges other than endpoints is a \emph{crossing}.
Based on the structure of $G^{4k}$, we have the following lemma.
\begin{lemma}\label{lemma43}
Let $\textrm{cr}(G^{4(k-1)})\geq 2(k-1)$ and let $D$ be a good drawing of $G^{4k}$ with $\textrm{cr}(D)<2k$.
In $D$, if there are $b$ $(b\geq 1)$ crossings on $P(v_i,v_{i+1},v_{i+2},v_{i+3},v_{i+4})$,
then $\textrm{cr}_D(v_iv_{i+4},G^{4k}\setminus E(P(v_i,v_{i+1},v_{i+2},v_{i+3},v_{i+4}))\geq b-1$, $i\in\mathbb{Z}_{4k}$.
\end{lemma}
\begin{proof}
Let $c=\textrm{cr}_D(v_iv_{i+4},G^{4k}\setminus E(P(v_i,v_{i+1},v_{i+2},v_{i+3},v_{i+4}))$.
Let $D_1$ be the drawing obtained by drawing a new edge $e$ joining $v_{i-4}$ and $v_{i+4}$ close enough to the path $v_{i-4}v_iv_{i+4}$ in $D$,
such that $e$ crosses an edge $e'$ if and only if $v_{i-4}v_iv_{i+4}$ crosses $e'$.
Let $D_2$ be the drawing obtained by deleting the edges in path $v_{i-4}v_iv_{i+1}v_{i+2}v_{i+3}v_{i+4}$ from $D_1$.
By the construction of $D_2$, it is easy to check that $D_2$ is a good drawing of a subdivision of $G^{4(k-1)}$
and $\textrm{cr}(D_2)=\textrm{cr}(D)-b+c$.
Then $c=b+\textrm{cr}(D_2)-\textrm{cr}(D)$.
Since subdivision does not affect the crossing number of a graph, we have $\textrm{cr}(D_2)\geq 2(k-1)$ by $\textrm{cr}(G^{4(k-1)})\geq 2(k-1)$.
By $\textrm{cr}(D)<2k$, we have $\textrm{cr}(D)\leq 2k-1$ and $-\textrm{cr}(D)\geq 1-2k$.
Then $c\geq b+2(k-1)+1-2k$. Hence $c\geq b-1$.
\end{proof}
\begin{remark}\label{re2}
The drawing in Fig. \ref{F1}b is obtained from the drawing $D$ in Fig. \ref{F1}a for $i=4$,
by the way to obtain $D_2$ in the proof of Lemma \ref{lemma43}.
\end{remark}
Let $D$ be a good drawing of graph $G$ and let $A$ be a subgraph of $G$.
Let $D(A)$ denote the sub-drawing of $A$ in $D$.
In $D$, an edge is \emph{clean} if it is not crossed by any other edges; otherwise, it is \emph{crossed}.
In $D$, $A$ is \emph{clean} if each edge in $A$ is clean.
For any good drawing $D$ of $G^{4k}$, we have $\textrm{cr}_D(H_i)=0$, $i\in\mathbb{Z}_{4k}$, since $D$ is a good drawing.
For $i\in\mathbb{Z}_{4k}$, by \eqref{fun}, we have
\begin{equation}\label{fun1}
f_D(H_i)=\textrm{cr}_D(H_i,G^{4k}\setminus E(H_i))/2.
\end{equation}
By \eqref{I1}, \eqref{I2} and \eqref{fun1}, we have the following lemma:
\begin{lemma}\label{lemma42}
Let $D$ be a good drawing of $G^{4k}$. Then $\textrm{cr}(D)=\sum_{i=0}^{4k-1}f_{D}(H_i)$.
\end{lemma}
Let $H_{i,j},i\leq j,$ be the union of $H_i,...,H_j$ with subscripts taken modulo $4k$,
and let $\textrm{cr}(G^{4(k-1)})\geq 2(k-1)$.
In the following of this section, assume that $D$ is a good drawing of $G^{4k}$ such that $\textrm{cr}(D)<2k$ and $f_D(H_{0,i})<(i+1)/2$ for all $i\in\mathbb{Z}_{4k}$.
The next two lemmas describe $\textrm{cr}_D(P(v_1,v_2,v_3,v_4),C(v_0,v_1,v_5,v_4))$.
\begin{lemma}\label{lemma411}
Let $\textrm{cr}(G^{4(k-1)})\geq 2(k-1)$ and let $D$ be a good drawing of $G^{4k}$ such that $\textrm{cr}(D)<2k$ and $f_D(H_{0,i})<(i+1)/2$ for all $i\in\mathbb{Z}_{4k}$.
Then $\textrm{cr}_D(P(v_1,v_2,v_3,v_4),C(v_0,v_1,v_5,v_4))\not= 0$.
\end{lemma}
\begin{lemma}\label{lemma48}
Let $\textrm{cr}(G^{4(k-1)})\geq 2(k-1)$ and let $D$ be a good drawing of $G^{4k}$ such that $\textrm{cr}(D)<2k$ and $f_D(H_{0,i})<(i+1)/2$ for all $i\in\mathbb{Z}_{4k}$.
Then $\textrm{cr}_D(P(v_1,v_2,v_3,v_4),C(v_0,v_1,v_5,v_4))\not=1$.
\end{lemma}
We postpone the proofs of Lemmas \ref{lemma411} and \ref{lemma48}.
With their proofs, we can prove the main result of this section.

\emph{Proof of Theorem \ref{T41}.}
By Lemma \ref{lemma41}, it suffices to prove that $\textrm{cr}(G^{4k})\geq 2k,k\geq 4$.
We prove it by induction on $k$.
By Lemma \ref{C16}, we have $\textrm{cr}(G^{4k})=2k,k=4$.
Let $k\geq 5$ and suppose that $\textrm{cr}(G^{4(k-1)})\geq 2(k-1)$.
To the contrary, suppose that $\textrm{cr}(G^{4k})< 2k$ and
$D$ is a good drawing of $G^{4k}$ with $\textrm{cr}(D)=\textrm{cr}(G^{4k})$.
Then $f_D(H_0),...,f_D(H_{4k-1})$ is a list of $4k$ real numbers,
and $\Sigma_{i=0}^{4k-1}f_D(H_i)=\textrm{cr}(D)$ by Lemma \ref{lemma42}.
By Theorem \ref{theo}, without loss of generality, assume that $f_D(H_{0,i})<(i+1)/2$, for all $i\in\mathbb{Z}_{4k}$.

Since $f_D(H_0)<1/2$, $v_0v_4$ is clean in $D$. Then by Lemma \ref{lemma43},
there is at most one crossing on $P(v_1,v_2,v_3,v_4)$ in $D$.
Hence $\textrm{cr}_D(P(v_1,v_2,v_3,v_4),C(v_0,v_1,v_5,$ $v_4))\leq 1$.
By Lemma \ref{lemma411}, $\textrm{cr}_D(P(v_1,v_2,v_3,v_4),C(v_0,v_1,v_5,v_4))=1$.
It contradicts Lemma \ref{lemma48}
$\hfill\Box$

By the Jordan Curve Theorem, we have the following lemma.
\begin{lemma}\label{J}
In a graph $G$, let $C$ and $C'$ be two vertex-disjoint cycles and $P_t=u_1u_2...u_t$ be a path of order $t$ with $V(P_t)\cap V(C)=\emptyset$.
Suppose that $D$ is a good drawing of $G$.
Then $\textrm{cr}_D(C,C')$ is even; $\textrm{cr}_D(C,P_t)$ is even if $u_1$ and $u_t$ are in the same region of $D(C)$, and it is odd otherwise.
\end{lemma}
The following result is used in both proofs of Lemmas \ref{lemma411} and \ref{lemma48}.
\begin{lemma}\label{lemma46}
Let $\textrm{cr}(G^{4(k-1)})\geq 2(k-1)$ and let $D$ be a good drawing of $G^{4k}$ such that $\textrm{cr}(D)<2k$ and $f_D(H_{0,i})<(i+1)/2$ for all $i\in\mathbb{Z}_{4k}$.
Then there are at most two crossings on $C(v_0,v_1,v_5,v_4)$ in $D$.
\end{lemma}
\begin{proof}
Since $f_D(H_0)<1/2$ and $f_D(H_{0,1})<1$, $v_1v_0v_4$ is clean and $v_1v_5$ is crossed at most once in $D$.
It suffices to prove that $v_4v_5$ is crossed at most once in $D$.
If there are at least two crossings on $v_4v_5$,
then each of $v_1v_5$, $v_2v_6$ and $v_3v_7$ is crossed by Lemma \ref{lemma43}, which contradicts $f_D(H_{0,4})<5/2$.
\end{proof}
The following two propositions show properties of $D$ if $\textrm{cr}_D(P(v_1,v_2,v_3,v_4),$ $C(v_0,v_1,v_5,v_4))=0$.
The first states that $D(P(v_1,v_2,v_3,v_4)\cup C(v_0,v_1,v_5,$ $v_4))$ is unique, up to isomorphism.
\begin{pro}\label{lemma49}
Let $\textrm{cr}(G^{4(k-1)})\geq 2(k-1)$ and let $D$ be a good drawing of $G^{4k}$ such that $\textrm{cr}(D)<2k$ and $f_D(H_{0,i})<(i+1)/2$ for all $i\in\mathbb{Z}_{4k}$.
If $\textrm{cr}_D(P(v_1,v_2,v_3,v_4),C(v_0,v_1,v_5,v_4))=0$,
then $D(P(v_1,v_2,v_3,v_4)\cup C(v_0,v_1,v_5,$ $v_4))$ is isomorphic to the drawing in Fig. \ref{F4}a.
\end{pro}
\begin{proof}
By $f_D(H_0)<1/2$, $v_1v_0v_4$ is clean in $D$.
Then $\textrm{cr}_D(C(v_0,v_1,v_5,v_4))=0$,
since $D$ is a good drawing.
If $\textrm{cr}_D(P(v_1,v_2,v_3,v_4))=0$, then it is trivial that $D(P(v_1,v_2,v_3,v_4)\cup C(v_0,v_1,v_5,v_4))$ is isomorphic to the drawing
in Fig. \ref{F4}a.
By contradiction. Suppose that $\textrm{cr}_D(P(v_1,v_2,v_3,v_4))>0$. Then $v_{1}v_{2}$ crosses $v_{3}v_{4}$, since $D$ is a good drawing.
Hence $D(P(v_1,v_2,v_3,v_4)\cup C(v_0,v_1,v_5,v_4))$ is isomorphic to the drawing in Fig. \ref{F4}b or Fig. \ref{F4}c.

In $D$, since $v_{1}v_{2}$ is crossed, $v_{1}v_5$ is clean by $f_D(H_{0,1})<1$,
and then $P(v_1,v_2,$ $v_3,v_4,v_5)$ cannot be crossed by any other edge by Lemma \ref{lemma43}.
In both cases, we obtain contradictions.
By Lemma \ref{J},
$v_2v_6v_5$ crosses $H_0\cup P(v_1,v_2,v_3,v_4)$ in Fig. \ref{F4}b,
and $v_0v_{4k-1}v_3$ crosses $C(v_1,v_2,v_3,v_4,v_5)$ in Fig. \ref{F4}c.
\end{proof}
The second states that in $D$, $C(v_1,v_2,v_3,v_4,v_5)$ crosses $C(v_{4k-5},v_{4k-4},v_0,$ $v_{4k-1})$ exactly twice
with $v_1v_5$ crossing $C(v_{4k-5},v_{4k-4},v_0,v_{4k-1})$ only once.
\begin{pro}\label{lemma410}
Let $\textrm{cr}(G^{4(k-1)})\geq 2(k-1)$ and let $D$ be a good drawing of $G^{4k}$ such that
$\textrm{cr}(D)<2k$ and $f_D(H_{0,i})<(i+1)/2$ for all $i\in\mathbb{Z}_{4k}$.
If $\textrm{cr}_D(P(v_1,v_2,v_3,v_4), C(v_0,v_1,v_5,v_4))=0$,
then $\textrm{cr}_D(C(v_1,v_2,v_3,v_4,v_5)
,C(v_{4k-5},$ $v_{4k-4},v_0,v_{4k-1}))=2$
and $\textrm{cr}_D(v_1v_5,C(v_{4k-5},v_{4k-4},v_0,v_{4k-1}))=1$.
\end{pro}
\begin{proof}
By Proposition \ref{lemma49}, the plane are divided into 3 regions $R_1, R_2$ and $R_3$ in $D(P(v_1,v_2,v_3,v_4)\cup C(v_0,v_1,v_5,v_4))$,
as shown in Fig. \ref{F4}a.
By $f_D(H_0)<1/2$, $v_1v_0v_4$ is clean.
Then by Lemma \ref{lemma43}, $P(v_1,v_2,v_3,v_4)$ is crossed at most once.
If $v_{4k-1}$ is in $R_2$, then both $v_{5}v_{9}v_{10}v_{14}...v_{4k-2}v_{4k-1}$ and $v_{5}v_{6}v_7v_{11}...v_{4k-1}$
cross $P(v_1,v_2,v_3,v_4)$ by Lemma \ref{J}, a contradiction.
If $v_{4k-1}$ is in $R_3$,
then each of $v_{3}v_{4k-1}, v_{2}v_{4k-2}v_{4k-1}$ and $v_{3}v_7...v_{4k-1}$ crosses $C(v_0,v_1,v_5,v_4)$,
which contradicts Lemma \ref{lemma46}.
Therefore, $v_{4k-1}$ is in $R_1$.
\begin{figure}
  \centering
	\subfigure[]{
    \label{fig:subfig:a}
    \includegraphics[width=0.23\textwidth]{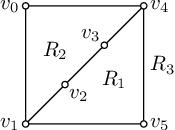}}
  \subfigure[]{
    \label{fig:subfig:a}
    \includegraphics[width=0.22\textwidth]{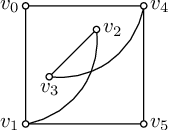}}
   \subfigure[]{
    \label{fig:subfig:a}
    \includegraphics[width=0.22\textwidth]{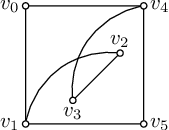}}
  \caption{$D(P(v_1,v_2,v_3,v_4)\cup C(v_0,v_1,v_5,$ $v_4))$.}\label{F4}
\end{figure}

By $f_D(H_{0,1})< 1$, $v_1v_5$ is crossed at most once in $D$.
By Lemma \ref{lemma43}, if $v_1v_5$ is crossed, then $P(v_1,v_2,v_3,v_4,v_5)$ is crossed at most twice in $D$;
if $v_1v_5$ is clean, then $P(v_1,v_2,v_3,v_4,v_5)$ is crossed at most once in $D$.
Hence $C(v_1,v_2,v_3,v_4,v_5)$ is crossed at most three times in $D$,
and then $\textrm{cr}_D(C(v_1,v_2,v_3,v_4,v_5),C(v_{4k-5},$ $v_{4k-4},v_0,v_{4k-1}))=2$ by Lemma \ref{J}.

Since $v_1v_5$ is crossed at most once, if $v_1v_5$ crosses $C(v_{4k-5},v_{4k-4},v_0,v_{4k-1})$,
then $\textrm{cr}_D(v_1v_5,C(v_{4k-5},v_{4k-4},v_0,v_{4k-1}))=1$.
By contradiction.
Suppose that $\textrm{cr}_D(v_1v_5,C(v_{4k-5},v_{4k-4},v_0,v_{4k-1}))=0$.
Then $\textrm{cr}_D(P(v_1,v_2,v_3,v_4,v_5),$ $C(v_{4k-5},v_{4k-4},v_0,v_{4k-1}))=2$.
Hence $v_1v_5$ is crossed by Lemma \ref{lemma43},
and then $v_1v_2$ is clean by $f_D(H_{0,1})<1$.
Since both $v_0v_4$ and $v_1v_2$ are clean, $v_2v_3v_4$ is crossed at most once by Lemma \ref{lemma43}.
Hence $v_4v_5$ crosses $C(v_{4k-5},v_{4k-4},v_0,v_{4k-1})$.
If $v_4v_5$ is crossed at least twice,
then both $v_2v_6$ and $v_3v_7$ are crossed by Lemma \ref{lemma43},
which contradicts $f_D(H_{0,4})<5/2$.
Thus both $v_4v_5$ and $v_2v_3v_4$ cross $C(v_{4k-5},$ $v_{4k-4},v_0,v_{4k-1})$.
If $v_3v_4$ is crossed, then both $v_2v_6$ and $v_3v_7$ are crossed by Lemma \ref{lemma43},
which contradicts $f_D(H_{0,3})< 2$.
If $v_2v_3$ is crossed, then $v_2v_6$ is crossed by Lemma \ref{lemma43},
which contradicts $f_D(H_{0,2})<3/2$.
\end{proof}
We are ready to give the proof of Lemma \ref{lemma411}.

\emph{Proof of Lemma \ref{lemma411}.}
By contradiction.
Suppose that $\textrm{cr}_D(P(v_1,v_2,v_3,v_4),$ $C(v_0,v_1,v_5,v_4))=0$.
By Proposition \ref{lemma410}, $\textrm{cr}_D(C(v_1,v_2,v_3,v_4,v_5),C(v_{4k-5},v_{4k-4},$ $v_0,v_{4k-1}))=2$
and $\textrm{cr}_D(v_1v_5,C(v_{4k-5},v_{4k-4},v_0,v_{4k-1}))=1$.
Thus $\textrm{cr}_D(C(v_{4k-5},$ $v_{4k-4},v_0,v_{4k-1}),C(v_1,v_2,v_6,v_5))\geq 2$ by Lemma \ref{J}.
Since $\textrm{cr}_D(v_1v_5,C(v_{4k-5},$ $v_{4k-4},v_0,v_{4k-1}))=1$,
$v_1v_2$ is clean and $v_1v_5$ cannot be crossed by any other edge in $D$ by $f_D(H_{0,1})< 1$.
Thus $v_2v_6v_5$ crosses $C(v_{4k-5},v_{4k-4},v_0,v_{4k-1})$ in $D$ by $\textrm{cr}_D(C(v_{4k-5},v_{4k-4},v_0,v_{4k-1}),C(v_1,v_2,v_6,v_5))\geq 2$.
Combining this result, we consider all remaining possible edges of $C(v_1,v_2,v_3,v_4,v_5)$ that crossing $C(v_{4k-5},v_{4k-4},v_0,v_{4k-1})$ in $D$.

Suppose that $v_{2}v_{3}$ crosses $C(v_{4k-5},v_{4k-4},v_0,v_{4k-1})$.
Then by $f_D(H_{0,2})<3/2$, $v_5v_6$ crosses $C(v_{4k-5},v_{4k-4},v_0,v_{4k-1})$.
Hence $v_2v_6$ is crossed by Lemma \ref{lemma43}, which contradicts $f_D(H_{0,2})<3/2$.

Suppose that $v_3v_4$ crosses $C(v_{4k-5},v_{4k-4},v_0,v_{4k-1})$.
If $v_5v_6$ is crossed, then both $v_2v_6$ and $v_3v_7$ are crossed by Lemma \ref{lemma43},
which contradicts $f_D(H_{0,3})< 2$.
Hence $v_5v_6$ is clean, and then $v_2v_6$ crosses $C(v_{4k-5},v_{4k-4},v_0,v_{4k-1})$.
By Lemma \ref{J}, we have $\textrm{cr}_D(C(v_2,v_3,v_7,v_6),C(v_{4k-5},v_{4k-4},v_0,v_{4k-1}))\geq 2$,
and hence $v_6v_7$ crosses $C(v_{4k-5},v_{4k-4},v_0,v_{4k-1})$ by $f_D(H_{0,3})< 2$.
Thus $v_3v_7$ is crossed by Lemma \ref{lemma43}, which contradicts $f_D(H_{0,3})<2$.

Suppose that $v_4v_5$ crosses $C(v_{4k-5},v_{4k-4},v_0,v_{4k-1})$.
Then by Lemma \ref{J}, $\textrm{cr}_D(C(v_3,v_4,v_5,v_6,v_7),C(v_{4k-5},v_{4k-4},v_0,v_{4k-1}))\geq 2$.
If $v_5v_6$ is crossed, then each of $v_2v_6$, $v_3v_7$ and $v_4v_8$ is crossed by Lemma \ref{lemma43},
which contradicts $f_D(H_{0,4})<5/2$.
Hence $v_5v_6$ is clean.
Similarly, $v_3v_4$ is clean and $v_4v_5$ cannot be crossed by any other edge in $D$.
Hence $v_3v_7v_6$ crosses $C(v_{4k-5},v_{4k-4},v_0,v_{4k-1})$
by $\textrm{cr}_D(C(v_3,v_4,v_5,v_6,v_7),C(v_{4k-5},v_{4k-4},v_0,v_{4k-1}))\geq 2$.
Since $v_5v_6$ is clean, $v_2v_6$ crosses $C(v_{4k-5},v_{4k-4},v_0,v_{4k-1})$.

If $v_{6}v_{7}$ is crossed, then both $v_{3}v_{7}$ and $v_{4}v_{8}$ are crossed by Lemma \ref{lemma43},
which contradicts $f_D(H_{0,4})<5/2$.
If $v_{3}v_{7}$ crosses $(C(v_{4k-5},v_{4k-4},v_0,v_{4k-1})$,
then $\textrm{cr}_D(C(v_{4k-5},v_{4k-4},v_0,v_{4k-1}),C(v_3,v_4,v_8,v_7))\geq 2$ by Lemma \ref{J}.
Thus $v_{7}v_8$ crosses $C(v_{4k-5},v_{4k-4},v_0,v_{4k-1})$ by $f_D(H_{0,4})<5/2$.
Then by Lemma \ref{lemma43}, $v_4v_8$ is crossed, which contradicts $f_D(H_{0,4})<5/2$.
$\hfill\Box$

\begin{figure}
  \centering
    \label{fig:subfig:a}
    \includegraphics[width=0.22\textwidth]{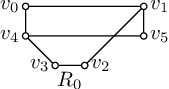}
  \caption{$D(P(v_1,v_2,v_3,v_4)\cup C(v_0,v_1,v_5,v_4))$.}\label{F3}
\end{figure}

The following proposition shows that $D(P(v_1,v_2,v_3,v_4)\cup C(v_0,v_1,v_5,v_4))$ is unique, up to isomorphism, if $\textrm{cr}_D(P(v_1,v_2,v_3,v_4),C(v_0,v_1,v_5,v_4))=1$.
\begin{pro}\label{lemma47}
Let $\textrm{cr}(G^{4(k-1)})\geq 2(k-1)$ and let $D$ be a good drawing of $G^{4k}$ such that
$\textrm{cr}(D)<2k$ and $f_D(H_{0,i})<(i+1)/2$ for all $i\in\mathbb{Z}_{4k}$.
If $\textrm{cr}_D(P(v_1,v_2,v_3,v_4),C(v_0,v_1,v_5,v_4))=1$,
then $D(P(v_1,v_2,v_3,v_4)\cup C(v_0,v_1,v_5,$ $v_4))$ is isomorphic to the drawing in Fig. \ref{F3}.
\end{pro}
\begin{proof}
Since $f_D(H_0)<1/2$, $v_1v_0v_4$ is clean.
Then $P(v_1,v_2,v_3,v_4)$ is crossed exactly once in $D$ by $\textrm{cr}_D(P(v_1,v_2,v_3,v_4),C(v_0,v_1,v_5,v_4))=1$ and Lemma \ref{lemma43},
and $\textrm{cr}_D(C(v_0,v_1,v_5,v_4))=0$ since $D$ is a good drawing.
We consider all possible edges of $P(v_1,v_2,v_3,v_4)$ that crossing $C(v_0,v_1,v_5,v_4)$ in $D$.

If $\textrm{cr}_D(v_2v_3,C(v_0,v_1,v_5,v_4))=1$, then $v_2$ and $v_3$ are in different regions of $D(C(v_0,v_1,v_5,v_4))$ by Lemma \ref{J}.
Hence both $v_2v_6v_7v_3$ and $v_2v_{4k-2}v_{4k-1}v_3$ cross $C(v_0,v_1,v_5,v_4)$ by Lemma \ref{J}.
Then there are at least three crossings on $C(v_0,v_1,v_5,v_4)$ in $D$,
which contradicts Lemma \ref{lemma46}.

If $\textrm{cr}_D(v_3v_4,C(v_0,v_1,v_5,v_4))=1$, then $\textrm{cr}_D(v_3v_4,v_1v_5)=1$,
since $D$ is a good drawing and $v_1v_0v_4$ is clean.
It follows that $\textrm{cr}_D(v_3v_4,C(v_1,v_2,v_6,v_5))=1$, since $P(v_1,v_2,v_3,v_4)$ is crossed exactly once in $D$.
Hence by Lemma \ref{J}, $v_3$ and $v_4$ are in different regions of $D(C(v_1,v_2,v_6,v_5))$.
Since $v_0v_4$ is clean, $v_3$ and $v_0$ are in different regions of $D(C(v_1,v_2,v_6,v_5))$.
By Lemma \ref{J}, both $v_3v_7v_8v_4$ and $v_0v_{4k-1}v_3$ cross $C(v_1,v_2,v_6,v_5)$.
Since $f_D(H_{0,1})<1$ and $v_1v_5$ is crossed,
both $v_1v_2$ and $v_1v_5$ cannot be crossed by any other edge in $D$.
Then both $v_3v_7v_8v_4$ and $v_0v_{4k-1}v_3$ cross $v_2v_6v_5$.
By $f_D(H_{0,2})<3/2$, at least one of them crosses $v_{5}v_{6}$.
Hence there are at least two crossings on $v_3v_4v_5v_6$.
By Lemma \ref{lemma43}, both $v_2v_{6}$ and $v_3v_7$ are crossed.
Then $f_D(H_{0,3})\geq 2$, a contradiction.

Therefore, $\textrm{cr}_D(v_1v_2,C(v_0,v_1,v_5,v_4))=1$.
Since $D$ is a good drawing and $v_0v_4$ is clean, $v_1v_2$ crosses $v_4v_5$.
Since $P(v_1,v_2,v_3,v_4)$ is crossed exactly once and $\textrm{cr}_D(C(v_0,v_1,v_5,v_4))=0$,
it is trivial that $D(P(v_1,v_2,v_3,v_4)\cup C(v_0,v_1,v_5,v_4))$ is
isomorphic to the drawing in Fig. \ref{F3}.
\end{proof}
The last lemma states that $f_D(H_{i,i+1})$ is at least $1/2$ for all $i\in\mathbb{Z}_{4k}$.
\begin{lemma}\label{lemma45}
Let $\textrm{cr}(G^{4(k-1)})\geq 2(k-1)$ and let $D$ be a good drawing of $G^{4k}$ such that
$\textrm{cr}(D)<2k$ and $f_D(H_{0,i})<(i+1)/2$ for all $i\in\mathbb{Z}_{4k}$.
Then $f_D(H_{i,i+1})\geq 1/2$, $i\in\mathbb{Z}_{4k}$.
\end{lemma}
\begin{proof}To the contrary, without loss of generality, suppose that $f_D(H_{0,1})=0$.
Then $v_1v_0v_4\cup v_2v_1v_5$ is clean.
Since $v_2v_1v_5$ is clean,
there is at most one crossing on $P(v_2,v_3,v_4,v_5)$ in $D$ by Lemma \ref{lemma43},
and then $\textrm{cr}_D(P(v_1,v_2,v_3,v_4)\cup C(v_0,v_1,v_5,v_4))\leq 1$.
If $\textrm{cr}_D(P(v_1,v_2,v_3,v_4)\cup C(v_0,v_1,v_5,v_4))=1$, then $v_2v_3$ crosses $v_4v_5$,
since $v_1v_0v_4\cup v_2v_1v_5$ is clean and $D$ is a good drawing.
Hence $\textrm{cr}_D(C(v_2,v_3,v_7,v_6),C(v_0,v_1,v_5,v_4))\geq 2$ by Lemma \ref{J}.
Then $v_4v_5$ crosses $C(v_2,v_3,v_7,v_6)$ at least twice, since $v_4v_0v_1v_5$ is clean.
It contradicts the fact that there is at most one crossing on $P(v_2,v_3,v_4,v_5)$.

Suppose that $\textrm{cr}_D(P(v_1,v_2,v_3,v_4)\cup C(v_0,v_1,v_5,v_4))=0$.
By Proposition \ref{lemma49}, $D(P(v_1,v_2,v_3,v_4)\cup C(v_0,v_1,v_5,v_4))$ is isomorphic to the drawing in Fig. \ref{F4}a.
We consider all possible regions $R_1$, $R_2$ and $R_3$ that $v_{4k-1}$ in, as shown in Fig. \ref{F4}a.
Since $v_1v_0v_4\cup v_2v_1v_5$ is clean, Lemma \ref{J} gives rise to the following results:
if $v_{4k-1}$ is in $R_1$,
then both $v_{4k-1}v_0$ and $v_{4k-1}v_{4k-5}v_{4k-4}v_{0}$ cross $P(v_2,v_3,v_4,v_5)$;
if $v_{4k-1}$ is in $R_2$,
then both $v_5v_{6}v_{10}...v_{4k-2}v_{4k-1}$ and $v_{5}v_{9}v_{13}...v_{4k-3}v_{4k-4}v_{4k-5}v_{4k-1}$
cross $v_2v_3v_4$;
if $v_{4k-1}$ is in $R_3$,
then both $v_{4k-1}v_3$ and $v_{4k-1}v_{4k-2}v_{2}$ cross $v_4v_5$.
In each case, there are at least two crossings on $P(v_2,v_3,v_4,v_5)$, a contradiciton.
\end{proof}

Now we are ready to give the proof of Lemma \ref{lemma48}.

\emph{Proof of Lemma \ref{lemma48}.}
To the contrary, suppose that $\textrm{cr}_D(P(v_1,v_2,v_3,v_4),$ $C(v_0,v_1,v_5,v_4))=1$.
By Proposition \ref{lemma47}, $D(P(v_1,v_2,v_3,v_4)\cup C(v_0,v_1,v_5,v_4))$ is isomorphic to the drawing in Fig. \ref{F3}.
$H_0$ is clean by $f_D(H_0)<1/2$.
Since $f_D(H_{0,1})<1$ and $v_{1}v_{2}$ is crossed, $H_{0,1}$ cannot be crossed by any other edge in $D$.
Thus $P(v_2,v_3,v_4,v_5)$ cannot be crossed by any other edge in $D$ by Lemma \ref{lemma43}.
Let $G_1=P(v_1,v_2,v_3,v_4)\cup C(v_0,v_1,v_5,v_4)$.
Since $G_1=H_{0,1}\cup P(v_2,v_3,v_4,v_5)$, $G_1$ cannot be crossed by any other edge in $D$.
If $v_{4k-1}$ is not in the region $R_0$ of $D(G_1)$ as shown Fig. \ref{F3},
then at least one of $v_0v_{4k-1}$ and $v_{4k-1}v_3$ crosses $G_1$ by Lemma \ref{J}, a contradiction.
Hence $v_{4k-1}$ is in $R_0$.

\begin{figure}
  \centering
    \subfigure[]{
    \label{fig:subfig:a}
    \includegraphics[width=0.22\textwidth]{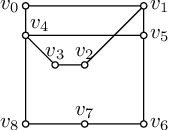}}
\subfigure[]{
    \label{fig:subfig:a}
    \includegraphics[width=0.22\textwidth]{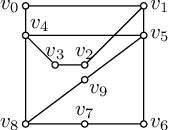}}
    \subfigure[]{
    \label{fig:subfig:a}
    \includegraphics[width=0.22\textwidth]{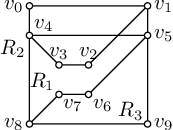}}
  \caption{Sub-drawings of $D$.}\label{F5}
\end{figure}

Let $\mathcal{P}_1$ be one of the two paths $v_2v_{4k-2}v_{4k-3}v_1$ and $v_3v_{4k-1}v_0$, and $\mathcal{P}_2$ be the other path.
Since $G_1$ cannot be crossed by any other edge in $D$, $\textrm{cr}_D(P(v_5,v_6,v_7,v_8)\cup v_4v_8,G_1)=0$,
and then both $\mathcal{P}_1$ and $\mathcal{P}_2$ cross $C(v_4,v_5,v_6,$ $v_7,v_8)$ by Lemma \ref{J}.
By Lemma \ref{lemma45}, $f_D(H_{2,3})\geq 1/2$.
Then by $f_D(H_{0,4})<5/2$, $v_4v_8$ is crossed at most once in $D$,
since $v_1v_2$ crosses $v_4v_5$.
Thus $P(v_4,v_5,v_6,v_7,v_8)$ is crossed at most twice in $D$ by Lemma \ref{lemma43}.
Since $v_4v_5$ cannot be crossed by any other edge in $D$, both $P(v_5,v_6,v_7,v_8)$ and $v_4v_8$ can be crossed at most once in $D$.
Thus, $\textrm{cr}_D(\mathcal{P}_1\cup \mathcal{P}_2,v_4v_8)=1$ and $\textrm{cr}_D(\mathcal{P}_1\cup \mathcal{P}_2,P(v_5,v_6,v_7,v_8))=1$.
Let $G_2=G_1\cup P(v_5,v_6,v_7,v_8)\cup v_4v_8$.
Then $G_2$ cannot be crossed by any other edge in $D$,
and $D(G_2)$ is isomorphic to the drawing in Fig. \ref{F5}a.

Let $G_3=G_2\cup v_5v_9v_8$.
Since $D$ is a good drawing and $G_2$ cannot be crossed by $v_5v_9v_8$,
$D(G_3)$ is isomorphic to the drawing in Fig. \ref{F5}b or \ref{F5}c.
In Fig. \ref{F5}b, both $v_2v_6$ and $v_3v_7$ are crossed by Lemma \ref{J}, which contradicts $f_D(H_{0,4})<5/2$.
Suppose that $D(G_3)$ is isomorphic to the drawing in Fig. \ref{F5}c and
$\mathcal{P}_1$ crosses $v_4v_8$.
Then by $f_D(H_{0,5})<3$, $v_5v_9$ is crossed at most once in $D$,
since $v_1v_2$ crosses $v_4v_5$ and $f_D(H_{2,3})\geq 1/2$.
Since $\mathcal{P}_2$ crosses $P(v_5,v_6,v_7,v_8)$, $v_8v_9$ is crossed at most once in $D$ by Lemma \ref{lemma43},
and $\textrm{cr}_D(\mathcal{P}_2,C(v_5,v_6,v_7,v_8,v_9))=2$ by Lemma \ref{J}.
Since $\mathcal{P}_2$ crosses $P(v_5,v_6,v_7,v_8)$ exactly once, $\mathcal{P}_2$ crosses $v_5v_9v_8$ exactly once.
Hence $v_5v_9v_8$ can be crossed at most once in $D$ except for $\mathcal{P}_2$.

Recall that
$\mathcal{P}_1$ is one of the two paths $v_2v_{4k-2}v_{4k-3}v_1$ and $v_3v_{4k-1}v_0$, and $\mathcal{P}_2$ is the other;
$G_3=G_2\cup v_5v_9v_8$, $\textrm{cr}_D(\mathcal{P}_1\cup \mathcal{P}_2,v_4v_8)=1$,
$\textrm{cr}_D(\mathcal{P}_1\cup \mathcal{P}_2,P(v_5,v_6,v_7,$ $v_8))=1$ and
$G_2$ cannot be crossed by any other edge in $D$ except for $\mathcal{P}_1\cup\mathcal{P}_2$;
$\textrm{cr}_D(\mathcal{P}_1\cup \mathcal{P}_2,v_5v_9v_8)=1$ and $v_5v_9v_8$ can be crossed at most once in $D$ except for $\mathcal{P}_1\cup\mathcal{P}_2$.
Hence $\textrm{cr}_D(\mathcal{P}_1\cup \mathcal{P}_2,C(v_5,v_6,v_7,v_8,v_9))=2$.

Moreover, recall that $v_{4k-1}$ is in $R_0$ of $D(G_1)$,
$v_{4k-1}$ is in one of the regions $R_1,R_2$ and $R_3$ of $D(G_3)$, as shown in Fig. \ref{F5}c.
If $v_{4k-1}$ is in $R_1$,
then $v_{4k-1}v_{4k-5}v_{4k-4}v_0$ crosses $C(v_4,v_5,v_6,v_7,v_8)$ by Lemma \ref{J},
a contradiction.
If $v_{4k-1}$ is in $R_2$,
then both $v_{6}v_{10}...v_{4k-2}v_{4k-1}$ and $v_{7}v_{11}...v_{4k-1}$ cross $v_5v_9v_8$ by Lemma \ref{J}, a contradiction.
Suppose that $v_{4k-1}$ is in $R_3$.
If $v_{4k-2}$ is in $R_3$, then $\textrm{cr}_D(\mathcal{P}_1\cup \mathcal{P}_2,C(v_5,v_6,v_7,v_8,v_9))\geq 4$ by Lemma \ref{J}, a contradiction.
If $v_{4k-2}$ is not in $R_3$,
then both $v_{4k-1}v_{4k-2}$ and $v_{4k-1}v_{4k-5}v_{4k-4}v_0$ cross $C(v_5,v_6,v_7,v_8,v_9)$, a contradiction.
$\hfill\Box$

\section{Acknowledgement.}

The work was supported by the National Natural Science
Foundation of China (No.\,61401186) and the Science
Foundation of the Educational Department of Liaoning Province (No.\,LJKZ0968).

\section{Reference}







\end{document}